\documentclass[11pt]{amsart}

\usepackage{amsmath,amsthm,amsfonts,amssymb,bm,graphicx}

\usepackage
{hyperref}
\hypersetup{colorlinks=true,citecolor=blue,linkcolor=blue,urlcolor=blue,
pdfstartview=FitH, pdfauthor=Valentin Blomer ,pdftitle=Applications of the Kuznetsov formula on $GL(3)$}

\theoremstyle{plain}
\newtheorem{theorem}{Theorem}
\newtheorem{cor}[theorem]{Corollary}
\newtheorem{lemma}{Lemma}
\newtheorem{prop}{Proposition}
\numberwithin{equation}{section}

\theoremstyle{definition}

\renewcommand{\geq}{\geqslant}
\renewcommand{\leq}{\leqslant}

\DeclareMathOperator{\SL}{SL}

\newcommand{\eps}{\varepsilon}
\newcommand{\RR}{\mathbb{R}}

\newcommand{\QQ}{\mathbb{Q}}

\newcommand{\NN}{\mathbb{N}}

\newcommand{\ov}[1]{\overline{#1}}

\makeatletter
\DeclareRobustCommand\widecheck[1]{{\mathpalette\@widecheck{#1}}}
\def\@widecheck#1#2{%
    \setbox\z@\hbox{\m@th$#1#2$}%
    \setbox\tw@\hbox{\m@th$#1%
       \widehat{%
          \vrule\@width\z@\@height\ht\z@
          \vrule\@height\z@\@width\wd\z@}$}%
    \dp\tw@-\ht\z@
    \@tempdima\ht\z@ \advance\@tempdima2\ht\tw@ \divide\@tempdima\thr@@
    \setbox\tw@\hbox{%
       \raise\@tempdima\hbox{\scalebox{1}[-1]{\lower\@tempdima\box
\tw@}}}%
    {\ooalign{\box\tw@ \cr \box\z@}}}
\makeatother

\begin{document}

\author{Valentin Blomer}
\address{Mathematisches Institut, Bunsenstr. 3-5, 37073 G\"ottingen, Germany} \email{blomer@uni-math.gwdg.de}

\title{Applications of the Kuznetsov formula on $GL(3)$}
 
\thanks{The author was supported   by a Volkswagen Lichtenberg Fellowship and a Starting Grant of the European Research Council.  }


\begin{abstract}  We develop a fairly explicit Kuznetsov formula on $GL(3)$ and discuss the analytic behaviour of the test functions on both sides.  Applications to Weyl's law, exceptional eigenvalues, a large sieve and $L$-functions are given. 
\end{abstract}


\setcounter{tocdepth}{2}  \maketitle 

\section{Introduction}

The   Bruggeman-Kuznetsov formula \cite{Br, Ku, DI} is one of the most powerful tools in the analytic theory of automorphic forms on $GL(2)$ and the cornerstone for the investigation of moments of families of $L$-functions, including striking applications to subconvexity and non-vanishing. It can be viewed as a relative trace formula for the group $G = GL(2)$ and the abelian subgroup $U_2 \times U_2 \subseteq G \times G $ where $U_2$ is the group of unipotent upper triangular matrices. The Kuznetsov formula in the simplest case is an equality of the shape
\begin{equation}\label{proto}
 2\sum_j \frac{\lambda_j(n) \overline{\lambda_j(m)}}{L({\rm Ad}^2 u_j, 1)}  h(t_j) +   \frac{1}{\pi} \int_{-\infty}^{\infty} \frac{\sigma_t(n)\overline{\sigma_t(m)}}{|\zeta(1+2it)|^2} h(t) dt =\delta_{n, m} \int_{-\infty}^{\infty} h(t)  d_{\text{spec}}t + \sum_{c} \frac{1}{c} S(n, m, c) h^{\pm} \left(\frac{|nm|}{c^2}\right)
\end{equation}
where
\begin{itemize}
\item $n, m \in \Bbb{Z} \setminus \{0\}$, 
\item $\delta_{n, m}$ is the Kronecker symbol,
\item the sum on the left-hand side runs over an orthogonal basis of Hecke-Maa{\ss} cusp forms $u_j$ for the group $SL_2(\Bbb{Z})$ having spectral parameter $t_j$ and Hecke eigenvalues $\lambda_j(n)$ for $n \in \Bbb{N}$  (and $\lambda_j(-n) := \pm \lambda_j(n)$ depending on whether $u_j$ is even or odd),
\item $\sigma_t(n)$ is the Fourier coefficient of an Eisenstein series defined by
\begin{displaymath}
  \sigma_t(n) = \sum_{d_1d_2 = |n|} d_1^{it} d_2^{-it}, 
\end{displaymath}
\item $d_{\text{spec}}t = \pi^{-2} t \tanh(\pi t)  d t$ is the spectral measure,
\item $h$ is some sufficiently nice, even test function, and
\item $h^{\pm}$ is a certain integral transform of $h$, the sign being $\text{sgn}(nm)$, described in \eqref{trafo}.  
\end{itemize}
There have been many generalizations  of the Kuznetsov to other groups of real rank one or products thereof, see e.g. \cite{MW, Re, BM}, the first of which covers also the groups $SL_2(\Bbb{C})$, $SO(n, 1)$ and $SU(2, 1)$; see also \cite{CLPSS, EGM} for interesting applications.   For the groups $GL(n)$, $n > 2$, Kuznetsov-type formulae are available \cite[Theorem 11.6.19]{Go}, \cite{Ye}, but they are in considerably less explicit form.  

The power of the $GL(2)$ Kuznetsov formula lies in the fact that one can choose arbitrary (reasonable) test functions on either side of the formula, and the relevant integral transforms are completely explicit in terms of Bessel functions. In fact, we have
\begin{equation}\label{trafo}
  h^{\pm}(x) = \int_0^{\infty} \mathcal{J}^{\pm}(t, x) h(t)  d_{\text{spec}} t
\end{equation}
where 
\begin{displaymath}
  \mathcal{J}^{\pm}(t, x) =2\pi i \sinh(\pi t)^{-1} \left\{\begin{array}{l} J_{2 i t}(4\pi \sqrt{x}) - J_{-2it}(4\pi \sqrt{x})\\ I_{2 i t}(4\pi \sqrt{x}) - I_{-2it}(4\pi \sqrt{x}) \end{array}\right\};
\end{displaymath}
this is best understood in terms of its Mellin transform:  
\begin{displaymath}
  \widehat{\mathcal{J}}^{\pm}(t, s) = \int_0^{\infty} \mathcal{J}^{\pm}(t, x) x^{s-1} dx =    \frac{G_t(s)}{G_t(1-s)} \mp \frac{G_t(1 + s)}{G_t(2 - s)} 
\end{displaymath}
where
\begin{displaymath}
  G_t(s) =  \pi^{-s} \Gamma\left(\frac{s+it}{2}\right) \Gamma\left(\frac{s-it}{2}\right). 
\end{displaymath}
In addition, this transform can be inverted and is essentially unitary:
\begin{displaymath}
  h(t) \approx \int_0^{\infty} \mathcal{J}^{\pm}(t, x) h^{\pm}(x) \frac{dx}{x}.
\end{displaymath}
There is a subtlety, as in the 
$+$ case the image of the map $h \mapsto h^{+}$ is not dense, but its complement is well-understood.  These formulas together with standard facts about Bessel functions make it possible to apply the Kuznetsov formula in both directions. Unfortunately, such explicit knowledge is not available for $GL(n)$, $n \geq 3$. \\

The aim of this article to provide a ``semi-explicit" version of the Kuznetsov formula for $GL(3)$ together with some careful analysis of the various terms occurring on both sides of the formula, and to give some applications in Theorems \ref{thm1}--\ref{thm4} below. On the way we will  prove a number of useful auxiliary results for $GL(3)$ Whittaker functions, Eisenstein series and Kloosterman sums that may be helpful for further investigation of $GL(3)$ automorphic forms.   
The proof of the Kuznetsov formula proceeds along classical lines: we compute the inner product of two Poincar\'e series in two ways: by spectral decomposition and by unfolding and computing the Fourier expansion of the Poincar\'e series. The latter has been worked out in great detail in \cite{BFG}. 

The \emph{spectral side}  \eqref{spectral} of the $GL(3)$ formula consists of three terms:
\begin{itemize}
\item the contribution of  the cuspidal spectrum,
\item the contribution of the minimal parabolic Eisenstein series,
\item the contribution of the maximal parabolic Eisenstein series.
\end{itemize}  
The \emph{arithmetic side} \eqref{arithmetic} contains four terms: 
\begin{itemize}
\item the diagonal contribution corresponding to the identity element in the Weyl group,
\item two somewhat degenerate terms\footnote{In all our applications they will be negligible.}  corresponding to $\left(\begin{smallmatrix} & & 1\\1 & & \\ & 1 &\end{smallmatrix}\right)$ and $\left(\begin{smallmatrix} & 1& \\ & &1 \\1 &  &\end{smallmatrix}\right)$,
\item the contribution of the long Weyl element. 
\end{itemize}
Interestingly, the two remaining elements in the Weyl group do not contribute as long as $n_1, n_2, m_1, m_2$ are non-zero; in fact, these two furnish the $GL(2)$ formula which is hidden in the degenerate terms of the Fourier expansion of the Poincar\'e series. On the arithmetic side, the various variables $n_1, n_2, m_1, m_2, D_1, D_2$ appearing in the test function are explicitly coupled similarly as on the right hand side  of \eqref{proto}. 

The spectral side \eqref{spectral} contains the weight function 
\begin{displaymath}
  h(\nu_1, \nu_2) = \left|\int_0^{\infty} \int_0^{\infty} \tilde{W}_{\nu_1, \nu_2}(y_1, y_2) F(y_1, y_2) \frac{dy_1\, dy_2}{(y_1y_2)^3}\right|^2
\end{displaymath}
where $\tilde{W}_{\nu_1, \nu_2}$ is a normalized Whittaker function on $GL(3)$ and $F$ is any  compactly supported test  function (or with sufficiently rapid decay at 0 and $\infty$ in both variables). In principle this integral transform is invertible: it has been shown in \cite{GK} that  a natural generalization of the Kontorovich-Lebedev transform inversion formula holds for Whittaker functions on $GL(n)$, hence we have a recipe to find a suitable $F$ to construct our favourite non-negative function $h$. Proceeding in this way would however considerably complicate the analysis of the arithmetic side, and hence we take a different route which is somewhat less precise, but more convenient for applications.  In Proposition \ref{whitint} below we show roughly the following: taking 
\begin{displaymath}
  F(y_1, y_2) =(\tau_1\tau_2(\tau_1+\tau_2))^{1/2}  f_1(y_1) f_2(y_2)   y_1^{ i(\tau_1+2\tau_2)} y_2^{ i(2\tau_1 + \tau_2)} 
\end{displaymath} 
for some fixed functions $f_1$, $f_2$ with compact support in $(0, \infty)$ 
yields a non-negative smooth bump function $h$ with $h(\nu_1, \nu_2) \asymp 1$ for $\nu_j = i\tau_j + O(1)$ and rapid decay outside this range. In other words, $h$ is a good approximation to the characteristic function of a unit square in the $(\nu_1, \nu_2)$-plane. Integration over $\tau_1, \tau_2$ can now give a good approximation to the characteristic function of any reasonable shape. Passing to a larger region in this way will in fact improve the performance of sum formula and ease the estimations on the arithmetic side. 

The test functions on the arithmetic side are completely explicit in \eqref{test1}, \eqref{test2} and  given as a multiple integral.  At least in principle a careful asymptotic analysis should yield a complete description of the behaviour of this function, but this seems very complicated. Nevertheless, we are able to give some non-trivial (and in some cases best possible) bounds in Proposition \ref{prop6} that suffice for a number of  applications that we proceed to describe.\\  

The commutative algebra $\mathcal{D}$ of invariant differential operators of $SL_3(\Bbb{R})$ acting on $L^2(SL_3(\Bbb{R})/SO_3)$ is generated by two elements (see \cite[p.\ 153]{Go}), the Laplacian and  another operator of degree 3. One class of  eigenfunctions of $\mathcal{D}$ is given by the power functions $I_{\nu_1, \nu_2}$ defined in \eqref{defI} below. A Maa{\ss} form $\phi$ for the group $SL_3(\Bbb{Z})$ with \emph{spectral parameters} $\nu_1, \nu_2$ is an element in $L^2(SL_3(\Bbb{Z})\backslash SL_3(\Bbb{R})/SO_3)$ that is an eigenfunction of $\mathcal{D}$ with the same eigenvalues as $I_{\nu_1, \nu_2}$ and vanishes along all parabolics, that is,
\begin{displaymath}
  \int_{(SL_3(\Bbb{Z}) \cup U)\backslash U} \phi(uz) du = 0
\end{displaymath} 
for $U = \left\{\left(\begin{smallmatrix} 1 &   & \ast\\  & 1 & \ast \\ & & 1 \end{smallmatrix}\right) \right\}$ and $\left\{\left(\begin{smallmatrix} 1 &  \ast & \ast\\  & 1 &  \\ & & 1 \end{smallmatrix}\right) \right\}$ (and then automatically for the minimal parabolic). We choose an orthonormal  basis $\{\phi_j\} = \mathcal{B} \subseteq L^2(SL_3(\Bbb{Z})\backslash SL_3(\Bbb{R})/SO_3)$ of Hecke-Maa{\ss} cusp forms (i.e.\ Maa{\ss} forms that are eigenfunctions of the Hecke algebra  as in  \cite[Section 6.4]{Go}) with spectral parameters $\nu^{(j)}_1, \nu_2^{(j)}$. If no confusion can arise, we drop the superscripts  $(j)$.

This can be re-phrased in more representation theoretic terms.  Let $SL_3(\Bbb{R}) = NAK$ be the Iwasawa decomposition where $K = SO_3$, $N$ is the standard unipotent subgroup and   
$A$ is the group of diagonal matrices with determinant 1 and positive entries, and let $\mathfrak{a}$ be the   Lie algebra of $A$. An infinite-dimensional, irreducible, everywhere unramified cuspidal automorphic representation $\pi$ of $GL_3(\Bbb{A}_{\Bbb{Q}})$ with trivial central character  is generated by a Hecke-Maa{\ss} form $\phi_j$ for $SL_3(\Bbb{Z})$ as above. The local (spherical) representation $\pi_{\infty}$ is an induced representation from the  parabolic subgroup $NA$ of the extension of a character  $\chi : A \rightarrow \Bbb{C}^{\times}$, $\text{diag}(x_1, x_2,x_3) \mapsto x_1^{\alpha_1} x_2^{\alpha_2} x_3^{\alpha_3}$ with $\alpha_1+\alpha_2+\alpha_3 = 0$. In this way we can identify the spherical cuspidal automorphic spectrum with a discrete subset of the Lie algebra $ \mathfrak{a}^{\ast}_{\Bbb{C}}/W$  ($W$ the Weyl group), where we associate to each Maa{\ss} form $\phi_j \in \mathcal{B}$ the linear form  $l = (\alpha_1, \alpha_2, \alpha_3) \in    \mathfrak{a}^{\ast}_{\Bbb{C}}/W$ that contains the (archimedean) Langlands parameters. A convenient basis in $\mathfrak{a}^{\ast}_{\Bbb{C}}$ is given by the fundmental weights $\text{diag}(2/3, -1/3, -1/3)$, $\text{diag}(1/3, 1/3, -2/3)$ of $SL_3$. The coefficients of $l= (\alpha_1, \alpha_2, \alpha_3)$ with respect to this basis can be obtained by  evaluating $l$ at the two co-roots $  
\text{diag}(1, -1, 0),  
 \text{diag}(0, 1, -1) \in \mathfrak{a}$ and are given by   $3\nu_1, 3\nu_2$.  
We then have $\alpha_1 = 2\nu_1+\nu_2$, $\alpha_2 = -\nu_1+\nu_2$, $\alpha_ 3 = -\nu_1 -2 \nu_2$. With this normalization, $\phi$ is an eigenform of the Laplacian with eigenvalue 
\begin{equation}\label{laplace}
\lambda = 1 - 3\nu_1^2 - 3\nu_1\nu_2 - 3\nu_2^2 = 1 - \frac{1}{2}(\alpha_1^2 + \alpha_2^2 + \alpha_3^2), 
\end{equation}
and the trivial representation is sitting at $(\nu_1, \nu_2) = (1/3, 1/3)$. 
The Ramanujan conjecture states that the Langlands parameters $\alpha_1, \alpha_2, \alpha_3$ of Maa{\ss} forms are purely  imaginary (equivalently, the spectral parameters $\nu_1, \nu_2$ are purely imaginary). A Maa{\ss} form is called exceptional if it violates the Ramanujan conjecture. Modulo the action of the Weyl group,  we can always assume that $\Im \nu_1, \Im \nu_2 \geq 0$ (positive Weyl chamber). Switching to the dual Maa{\ss} form if necessary, we can even assume without loss of generality $0 \leq \Im \nu_1 \leq \Im \nu_2$. \\

A count of the Maa{\ss} forms $\phi \in \mathcal{B}$  inside the ellipse $\lambda \leq T^2$ described by \eqref{laplace} is referred to as Weyl's law. The number of such forms is known to be  $cT^{5} + O(T^3)$ for some constant $c$, see  
\cite{Mi, LM}.  As a first test case of the Kuznetsov formula we show a result of comparable strength as \cite[Proposition 4.5]{LM} that turns out to be a simple corollary of the Kuznetsov formula. A similar upper bound   has recently been proved   by X.\ Li \cite{Li}.   Let $L(\phi \times \tilde{\phi}, s)$ be the  Rankin-Selberg $L$-function (see \eqref{RS} below). Then the following  weighted count of the cuspidal spectrum in a small ball of radius $O(1)$ in $\mathfrak{a}_{\Bbb{C}}^{\ast}$ holds. 

\begin{theorem}\label{thm1} There are absolute constants $ c_1, c_2> 0$, $T_0, K \geq 1$ with the following property:  for all  $T_1, T_2  \geq T_0$  we have 
\begin{displaymath}
c_1 T_1T_2(T_1+T_2)\leq   \sum_{\substack{| \nu^{(j)}_1 -iT_1| \leq K\\  | \nu^{(j)}_2 - iT_2| \leq K}} \left(\underset{s=1}{\rm res} L(\phi_j \times \tilde{\phi}_j, s)\right)^{-1}  \leq c_2 T_1T_2(T_1+T_2). 
  \end{displaymath}  
\end{theorem}

It is standard to estimate the residue from above, but due to possible Siegel zeros a good lower bound is not known. If $\phi = \text{sym}^2 u$ for some Hecke-Maa{\ss} form $u \in L^2(SL_2(\Bbb{Z})\backslash \mathfrak{h}^2)$ with spectral parameter $\nu \in i\Bbb{R}$, then Ramakrishnan and Wong \cite{RW} have shown that no Siegel zeros exist:
\begin{displaymath}
  \underset{s=1}{\rm res} L(\text{sym}^2u \times \text{sym}^2u, s)  = (1+|\nu|)^{o(1)}. 
\end{displaymath}
In general we will only be able to prove the following bounds: if $\phi$ has spectral parameters $\nu_1, \nu_2$, then setting   $C := (1+|\nu_1 + \nu_2|) (1+|\nu_1|) (1+|\nu_2|)$ we have 
\begin{equation}\label{resiweak}
C^{-1}\ll 
  \underset{s=1}{\rm res} L(\phi \times \tilde{\phi}, s) \ll C^{\varepsilon}.
 \end{equation}
In particular it follows (after possibly enlarging the constant $K$ in Theorem \ref{thm1}) that  in each ball inside $i \mathfrak{a}^{\ast}$ of sufficiently large constant radius, there exist cusp forms. We will prove \eqref{resiweak} in Lemma \ref{lem2} below. \\

Miller \cite{Mi} proved that almost all forms are non-exceptional, that is, the number  of exceptional forms $\phi_j \in \mathcal{B}$ with $\lambda_j \leq T^2$  is $o(T^5)$. This was, among other things, strengthened in \cite{LM} to $O(T^3)$.  By unitaricity and the standard Jacquet-Shalika bounds towards the Ramanujan conjecture\footnote{better bounds are available by the work of Luo-Rudnick-Sarnak, but this is not needed here. Even the value of the constant $1/2$ is irrelevant.}   (cf.\ \eqref{alpha1} below)  the spectral parameters $\nu_1, \nu_2$ of an exceptional Maa{\ss} form are of the form (assuming  $0 \leq \Im \nu_1 \leq \Im\nu_2$)
\begin{displaymath}
  (\nu_1, \nu_2) = (2\rho/3, -\rho/3 + i\gamma), \quad \gamma \geq 0, \, |\rho| \leq 1/2,
\end{displaymath}
see \eqref{nu1} below. It is an easy corollary of Theorem \ref{thm1} that there are  $O(T^{2+\varepsilon})$ exceptional eigenvalues with $\gamma = T + O(1)$, but   more can be shown which can be viewed as a density theorem for exceptional eigenvalues and   interpolates nicely between the Jacquet-Shalika bounds and the tempered spectrum.

\begin{theorem}\label{thm2} For any $\varepsilon > 0$ we have
\begin{displaymath}
   \sum_{\substack{ \phi_j \text{ exceptional}\\ \gamma_j = T+O(1)}}  T^{4|\rho_j|} \ll_{\varepsilon} T^{2+\varepsilon}. 
   \end{displaymath}
      \end{theorem}
 
Next we prove a large sieve type estimate for Hecke eigenvalues. Let $A_j(n, 1)$ denote the Hecke eigenvalues of the Hecke-Maa{\ss} cusp form $\phi_j$. 

\begin{theorem}\label{thm3} Let $N \geq 1$, $T_1, T_2 \geq T_0$ sufficiently large, and let $\alpha(n)$ be a   sequence of complex numbers. Then
\begin{equation}\label{large}
 \sum_{\substack{T_1 \leq |\nu^{(j)}_1| \leq 2T_1\\ T_2 \leq  |\nu^{(j)}_2| \leq 2 T_2}} \Bigl|\sum_{n \leq N} \alpha(n) A_j(n, 1)\Bigr|^2 \ll_{\varepsilon} (T_1^2T_2^2(T_1+T_2) + T_1T_2N^2)^{1+\varepsilon} \|\alpha \|_2^2
\end{equation}
for any $\varepsilon> 0$ where $\|\alpha \|_2 = (\sum_n |\alpha(n)|^2)^{1/2}$. 
\end{theorem}

The first term is optimal on the right hand side is optimal.  Most optimistically one could hope for an additional term of size $N$ (instead of $T_1T_2N^2$), but in any case our result suffices for   an essentially  optimal  bound of the second moment of a  family of genuine $GL(3)$ $L$-functions. This seems to be the first bound of this kind in the literature. For large sieve inequalities in the level aspect (with very different proofs) see \cite[Theorem 4]{DK} and \cite{Ve}. 

\begin{theorem}\label{thm4} For $T  \geq 1$ and any $\varepsilon > 0$ we have 
\begin{displaymath}
   \sum_{ T \leq |\nu^{(j)}_1|,   |\nu^{(j)}_2| \leq  2T} |L(\phi_j, 1/2)|^2 \ll_{\varepsilon}  T^{5+\varepsilon}.
\end{displaymath}
\end{theorem}

More applications of the $GL(3)$ Kuznetsov formula to the Sato-Tate distribution of $GL(3)$ Hecke eigenvalues and a version of Theorem \ref{thm2} for the Langlands parameters at finite places will be given in a forthcoming paper \cite[Theorems 1 - 3]{BBR}. \\

After the paper was submitted, two other    interesting approaches to the $GL(3)$ Kuznetsov formula have been developed independently by Buttcane \cite{But} and Goldfeld-Kontorovich \cite{GK2}. The present technique, however,  gives   the strongest bounds for the Kloosterman terms in the Kuznetsov formula which are indispensable for applications to $L$-functions as in Theorems \ref{thm3} and \ref{thm4}. One may compare, for instance, with \cite{GK2} for which the reader is referred to the appendix which features in Theorem \ref{appendix} another result of independent interest.\\

It would be very interesting to generalize the present results and techniques to congruence subgroups of $SL_3(\Bbb{Z})$ of the type 
$$\Gamma_0(q) = \left\{ \gamma \in \SL_3(\Bbb{Z}) \mid \gamma \equiv  \left(\begin{matrix} \ast & \ast & \ast \\ \ast & \ast & \ast \\  0 & 0 & \ast\end{matrix} \right) \, \, (\text{mod } q)\right\}.$$
The analytic parts of the present argument (in particular the bounds for Whittaker functions and the corresponding integral transforms) work without any change. One needs a more general Bruhat decomposition to calculate the Fourier expansion of the relevant Poincar\'e series, and it would be useful to have an explicit spectral decomposition for the space $L^2(\Gamma_0(q)\backslash  \mathfrak{h}^3)$. This along with further applications will be addressed in \cite{Bal}. \\

\textbf{Acknowledgement.} I would like to thank F.\ Brumley, J.\ Buttcane, A.\ Kontorovich, and F.\ Shahidi for enlightening discussions and for answering my questions on various aspects of automorphic forms and pointing out errors in an earlier version. Particularly, I would like to thank P.\ Sarnak for suggesting the application given in Theorem \ref{thm2}. The importance of deriving a Kuznetsov formula for $GL(3)$ that is user-friendly for analytic number theorists has been discussed at the workshop ``Analytic theory of $GL(3)$ automorphic forms and applications" at the American Institute of Mathematics in November 2008. In particular,  the questions of finding bounds in the situation of  Theorem \ref{thm3} and  Lemma \ref{lem2} have been listed as open problems. Finally, thanks are due to referee for useful comments.

\section{Whittaker functions}

Let $\nu_1, \nu_2 \in \Bbb{C}$. We introduce the notation
\begin{equation}\label{defnu1}
   \nu_0 := \nu_1 + \nu_2
 \end{equation}
and (as in the introduction)
\begin{equation}\label{defnu2}
  \alpha_1 = 2\nu_1 +\nu_2, \quad \alpha_2 = -\nu_1 + \nu_2, \quad \alpha_3 = -\nu_1 - 2\nu_2.
  \end{equation}
 The transformations
\begin{equation}\label{invariant}
    (\nu_1, \nu_2) \rightarrow (-\nu_1, \nu_0) \rightarrow (\nu_2, -\nu_0) \rightarrow (-\nu_2, -\nu_1) \rightarrow (-\nu_0, \nu_1) \rightarrow (\nu_0, -\nu_2) 
\end{equation}
leave $\{\alpha_1, \alpha_2, \alpha_3\}$ invariant, and they also leave $\{|\Im \nu_0|, |\Im \nu_1|, |\Im \nu_2|\}$ invariant. For convenience we   assume  the Jacquet-Shalika  bounds towards the Ramanujan conjecture 
\begin{equation}\label{alpha1}
  \max(|\Re \alpha_1|,  |\Re \alpha_2|, | \Re \alpha_3|)  \leq    \frac{1}{2}, 
\end{equation}  
and we always assume unitaricity
\begin{equation}\label{alpha2}
  \{\alpha_1, \alpha_2, \alpha_3\} = \{-\overline{\alpha_1}, - \overline{\alpha_2}, - \overline{\alpha_3}\}.
\end{equation}
It is elementary to deduce from \eqref{alpha1} that
\begin{equation}\label{nu0}
  \max(|\Re \nu_0|,  |\Re \nu_1|, | \Re \nu_2|)  \leq \frac{1}{3} 
\end{equation}
and to deduce from \eqref{alpha2} that
\begin{equation}\label{nu}
  \nu_0, \nu_1, \nu_2, \alpha_1, \alpha_2, \alpha_3 \in i\Bbb{R}
  \end{equation}
 or
 \begin{equation}\label{nu1} 
 \begin{split}
 &  \{\alpha_1, \alpha_2, \alpha_3\} \in \{\rho+i\gamma, -\rho+i\gamma, -2i\gamma  \}, \\
&  \{\nu_1, \nu_2, \nu_0\} \in \{2\rho/3,   - \rho/3 + i\gamma, \rho/3 + i \gamma\} \text{ or its translates under \eqref{invariant}}
   \end{split}
\end{equation}
with $\rho, \gamma \in \Bbb{R}$ and $|\rho| \leq 1/2$ by \eqref{alpha1}. The choice
\begin{equation}\label{nontemp}
 (\alpha_1, \alpha_2, \alpha_3) =  (\rho+i\gamma, -\rho+i\gamma,  -2 i\gamma), \quad  (\nu_1, \nu_2, \nu_0) = (2\rho/3,   - \rho/3 + i\gamma, \rho/3 + i \gamma)
\end{equation}
is unique if we require $\Im \nu_2 \geq \Im \nu_1 \geq 0$, $\gamma \geq 0$.   \\

 Let  $\mathfrak{h}^2 =  \left\{z= \left(\begin{matrix} 1 & x  \\ & 1\end{matrix} \right)  \left(\begin{matrix} y  &    \\ &  1\end{matrix} \right)  \mid y  > 0, \, x  \in \Bbb{R}\right\} \cong GL_2(\Bbb{R})/(O_2 Z_2) \cong SL_2(\Bbb{R})/SO_2$  and 
\begin{displaymath}
  \mathfrak{h}^3 = \left\{z= \left(\begin{matrix} 1 & x_2 & x_3\\ & 1 & x_1\\ & & 1\end{matrix} \right)  \left(\begin{matrix} y_1y_2 &  &  \\ & y_1 & \\ & & 1\end{matrix} \right) \mid y_1, y_2 > 0, \, x_1, x_2, x_3 \in \Bbb{R}\right\}  \cong GL_3(\Bbb{R})/(O_3Z_3) \cong SL_3(\Bbb{R})/SO_3. 
\end{displaymath}  
The group $SL_3(\Bbb{Z})$ acts faithfully on $\mathfrak{h}^3$ by left multiplication.

 The  Whittaker function $\mathcal{W}^{\pm}_{\nu_1, \nu_2} : \mathfrak{h}^3 \rightarrow \Bbb{C}$ is given by\footnote{Some authors use different signs in the long Weyl element, but since the $I_{\nu_1, \nu_2}$ function depends only on $y_1, y_2$, this leads to the same definition.} (analytic continuation in $\nu_1, \nu_2$ of) 
\begin{equation}\label{generalWhit}
  \mathcal{W}^{\pm}_{\nu_1, \nu_2}(z) = \int_{\Bbb{R}^3} I_{\nu_1, \nu_2}\left(\left(\begin{matrix}  & &1\\ & 1& \\  1& & \end{matrix}\right) \left(\begin{matrix} 1 & u_2& u_3\\ & 1 & u_1 \\  & & 1\end{matrix}\right)z\right) e(-u_1 \mp  u_2) du_1 \, du_2 \, du_3
 \end{equation}
with 
\begin{equation}\label{defI}
  I_{\nu_1, \nu_2}(z) = y_1^{1 + 2\nu_1 + \nu_2} y_2^{1+\nu_1+2\nu_2}.
\end{equation}  
   Compared to   \cite{Bu}\footnote{In \cite[p.\ 154, third display]{Go} the values of $\nu_1, \nu_2$ are interchanged in the definition of $I$-function, but the following formulas are again in accordance with Bump's definition.} we have re-normalized the indices $\nu_j \rightarrow 1/3 +  \nu_j$.    By the formula of Takhtadzhyan-Vinogradov  we have $ \mathcal{W}^{\pm}_{\nu_1, \nu_2}(z) = e(x_1 \pm x_2) W_{\nu_1, \nu_2}(y_1, y_2)$ where\footnote{The normalization is complicated: the leading constant in \cite[(6.1.3)]{Go} should be 8 instead of 4, while the definition \cite[(1.1)]{St} differs from \eqref{generalWhit}, in addition to the Gamma-factors,  by a factor 2.} 
   \begin{equation}\label{defW}
\begin{split}
  W_{\nu_1, \nu_2}(y_1, y_2) = & 8   y_1y_2 \left(\frac{y_1}{y_2}\right)^{\frac{\nu_1 - \nu_2}{2}}   \prod_{j=0}^2  \pi^{\frac{1}{2} + \frac{3}{2} \nu_j} \Gamma\left(\frac{1}{2} + \frac{3}{2} \nu_j\right)^{-1}  \\
  &\times \int_0^{\infty} K_{\frac{3}{2}\nu_0}(2\pi y_2\sqrt{1+1/u^2})  K_{\frac{3}{2}\nu_0}(2\pi y_1\sqrt{1+u^2}) u^{\frac{3}{2}(\nu_1-\nu_2)} \frac{du}{u}.
  \end{split}
\end{equation}
It is convenient to slightly re-normalize this function: let
\begin{equation}\label{const}
  c_{\nu_1, \nu_2} := \pi^{-3\nu_0} \prod_{j=0}^2    \Gamma\left(\frac{1}{2} + \frac{3}{2}  \nu_j\right) \left|\Gamma\left(\frac{1}{2} + \frac{3}{2} i \Im \nu_j\right)\right|^{-1}
\end{equation}
and
\begin{displaymath}
\begin{split}
   \tilde{W}_{\nu_1, \nu_2}(y_1, y_2) := & W_{\nu_1, \nu_2}(y_1, y_2)  c_{\nu_1, \nu_2}
  =   8\pi^{\frac{3}{2} }    \prod_{j=0}^2  \left|\Gamma\left(\frac{1}{2} + \frac{3}{2} i \Im \nu_j\right)\right|^{-1}   y_1y_2 \left(\frac{y_1}{y_2}\right)^{\frac{\nu_1 - \nu_2}{2}} \\
  &\times \int_0^{\infty} K_{\frac{3}{2}\nu_0}(2\pi y_2\sqrt{1+1/u})  K_{\frac{3}{2}\nu_0}(2\pi y_1\sqrt{1+u}) u^{\frac{3}{4}(\nu_1-\nu_2)} \frac{du}{u}.
  \end{split}
\end{displaymath}
If $\nu_1, \nu_2 \in i \Bbb{R}$, this changes the original Whittaker function only by a constant on the unit circle, in the situation \eqref{nontemp} it changes the order of magnitude by a bounded factor. Often the Whittaker function is defined entirely without the normalizing Gamma-factors in the denominator of \eqref{defW} in which case it is often referred to as the completed Whittaker function.  It is convenient not to work with the completed Whittaker function  here, see Remark 3 below. (Of course, $\tilde{W}_{\nu_1, \nu_2}$ is not analytic in the indices any more.)

The $GL(2)$-analogue of this function is
\begin{equation}\label{gl2}
 W_{\nu}(y) := 2 \pi^{\frac{1}{2} } \left|\Gamma\left(\frac{1}{2} + \nu\right)\right|^{-1} \sqrt{y} K_{\nu}(2 \pi y),
\end{equation} 
see \cite[p.\ 65]{Go}. 

We proceed to collect analytic information on the $GL(3)$ Whittaker function. 
We have the double Mellin inversion formula (\cite[p.\ 155]{Go}, \cite[(10.1)]{Bu})
\begin{equation}\label{doubleMellin}
\begin{split}
  \tilde{W}_{\nu_1, \nu_2}(y_1, y_2) & = \frac{y_1y_2 \pi^{\frac{3}{2}  } }{  \prod_{j=0}^2  |\Gamma\left(\frac{1}{2} + \frac{3}{2}i \Im \nu_j\right)|}\\
  & \times \frac{1}{(2\pi i)^2} \int_{(c_2)} \int_{(c_1)} \frac{\prod_{j=1}^3 \Gamma(\frac{1}{2}(s_1 + \alpha_j)) \prod_{j=1}^3 \Gamma(\frac{1}{2}(s_2 - \alpha_j)) }{4\pi^{s_1+s_2} \Gamma(\frac{1}{2}(s_1+s_2))} y_1^{-s_1} y_2^{-s_2} ds_1\, ds_2. 
  \end{split}
\end{equation}
This implies in particular that $\tilde{W}_{\nu_1, \nu_2}$ is invariant under the transformations \eqref{invariant} which is the reason for including the normalization constant $c_{\nu_1, \nu_2}$.  Note that
\begin{equation}\label{bar}
  \tilde{W}_{\nu_1, \nu_2}(y_1, y_2) = \tilde{W}_{\nu_2, \nu_1}(y_2, y_1) = \overline{\tilde{W}_{\overline{\nu_1}, \overline{\nu_2}}(y_1, y_2)}. 
\end{equation}

Uniform bounds for Bessel functions are rare in the literature, but frequently needed in the $GL(2)$ theory. We are not aware of any uniform bound for a $GL(n)$ Whittaker function with $n > 2$. Although the proofs of Theorems \ref{thm1} - \ref{appendix} do not require bounds for individual Whittaker functions, we record here for future reference   the following uniform result.  

\begin{prop}\label{prop1} Let $\nu_1, \nu_2 \in \Bbb{C}$ satisfy \eqref{nu0} - \eqref{nu1}   and write  $\theta = \max(|\Re \alpha_1|, |\Re \alpha_2|, |\Re \alpha_3|) \leq 1/2$.  Let $\theta  < \sigma_1 < \sigma_2$ and $\varepsilon > 0$.  Then for  any $\sigma_1 \leq c_1, c_2 \leq \sigma_2 $ we have 
\begin{displaymath}
  \tilde{W}_{\nu_1, \nu_2}(y_1, y_2) \ll_{\sigma_1, \sigma_2, \varepsilon} \frac{y_1y_2}{(1+|\nu_1| + |\nu_2|)^{1/2-\varepsilon}}\left(\frac{y_1}{ 1+|\nu_1| + |\nu_2|}\right)^{-c_1}  \left(\frac{y_2}{ 1+|\nu_1| + |\nu_2|}\right)^{-c_2}.
\end{displaymath}
 \end{prop}
\textbf{Remark 1.} This result can be refined somewhat, in particular for small and large $y_1, y_2$. In the ``transitional range" it is not too far from the truth. For instance, if $\nu_1 = \nu_2 = iT$ are large (and purely imaginary) and $y_1 = y_2 = \frac{3}{2\pi}T - \frac{1}{100} T^{1/3}$, then the integral in \eqref{defW} is non-oscillating, and it follows from  \eqref{defW} and known properties of the $K$-Bessel function that $\tilde{W}_{\nu_1, \nu_2}(y_1, y_2) \asymp y_1y_2/T^{4/3}$, whereas our bound gives $\tilde{W}_{\nu_1, \nu_2}(y_1, y_2) \ll y_1y_2T^{\varepsilon - \frac{1}{2}}$. \\

\textbf{Proof.}  
Let us first assume that $\nu_1, \nu_2$ are purely imaginary; we write $\tilde{\nu}_j := \Im \nu_j$.  By \eqref{bar} and the invariance under  \eqref{invariant}   we can assume without loss of generality that  \begin{equation}\label{range}
0 \leq   \tilde{\nu}_1 \leq \tilde{\nu}_2.
\end{equation}
 By \eqref{doubleMellin} and Stirling's formula, we have
\begin{equation}\label{stir}
\begin{split}
 \tilde{W}_{\nu_1, \nu_2}&(y_1, y_2)  \ll_{\sigma_1, \sigma_2} y_1y_2 \int_{-\infty}^{\infty}\int_{-\infty}^{\infty}  \frac{\prod_{j=1}^3 (1+ |it_1 + \alpha_j|)^{\frac{c_1-1}{2}} \prod_{j=1}^3 (1+|i t_2 - \alpha_j|)^{\frac{c_2-1}{2} }}{ (1+|t_1+t_2|)^{\frac{c_1+c_2-1}{2}} }  \\
  & \times \exp\left(\frac{3\pi}{4}\sum_{j=0}^2|\nu_j|   - \frac{\pi}{4}\sum_{j=1}^3|it_1 + \alpha_j| - \frac{\pi}{4}\sum_{j=1}^3|it_2 - \alpha_j| + \frac{\pi}{4} |t_1+t_2|\right) y_1^{-c_1} y_2^{-c_2} dt_1\, dt_2.
  \end{split}
\end{equation} 
for $\sigma_1 < c_1, c_2 < \sigma_2$. 
It is elementary to check that 
\begin{equation}\label{exp}
  \frac{3\pi}{4}\sum_{j=0}^2|\nu_j|   - \frac{\pi}{4}\sum_{j=1}^3|it_1 + \alpha_j| - \frac{\pi}{4}\sum_{j=1}^3|it_2 - \alpha_j| + \frac{\pi}{4} |t_1+t_2| \leq 0
\end{equation}
with equality if and only if
\begin{equation}\label{case1}
  \tilde{\nu}_1-\tilde{\nu}_2 \leq t_1 \leq \tilde{\nu}_1 + 2\tilde{\nu}_2, \quad \tilde{\nu}_2 -\tilde{\nu}_1 \leq t_2 \leq 2\tilde{\nu}_1 + \tilde{\nu}_2,
\end{equation}
or
\begin{equation}\label{case2}
  -2\tilde{\nu}_1-\tilde{\nu}_2 \leq t_1 \leq \tilde{\nu}_1 - \tilde{\nu}_2, \quad -\tilde{\nu}_1 -2\tilde{\nu}_2 \leq t_2 \leq \tilde{\nu}_2 - \tilde{\nu}_1.
\end{equation}
For $a < b$ let $w_{a, b}(t)$ be defined by
\begin{displaymath}
  w_{a, b} :=  \min(1,  e^{-\frac{\pi}{4}(t-b)}, e^{-\frac{\pi}{4}(a-t)});
\end{displaymath}
then the exp-factor in \eqref{stir} is bounded by 
\begin{displaymath}
  w_{  \tilde{\nu}_1-\tilde{\nu}_2,  \tilde{\nu}_1 + 2\tilde{\nu}_2}(t_1) w_{   \tilde{\nu}_2 -\tilde{\nu}_1, 2\tilde{\nu}_1 + \tilde{\nu}_2}(t_2) + 
   w_{ -2\tilde{\nu}_1-\tilde{\nu}_2, \tilde{\nu}_1 - \tilde{\nu}_2}(t_1) w_{ -\tilde{\nu}_1 -2\tilde{\nu}_2 , \tilde{\nu}_2 - \tilde{\nu}_1}(t_2),
\end{displaymath}
and hence we have
\begin{displaymath}
\begin{split}
   \tilde{W}_{\nu_1, \nu_2}(y_1, y_2) & \ll y_1y_2 \int_{-\infty}^{\infty}\int_{-\infty}^{\infty}  \left(w_{  \tilde{\nu}_1-\tilde{\nu}_2,  \tilde{\nu}_1 + 2\tilde{\nu}_2}(t_1) w_{   \tilde{\nu}_2 -\tilde{\nu}_1, 2\tilde{\nu}_1 + \tilde{\nu}_2}(t_2) + 
   w_{ -2\tilde{\nu}_1-\tilde{\nu}_2, \tilde{\nu}_1 - \tilde{\nu}_2}(t_1) w_{ -\tilde{\nu}_1 -2\tilde{\nu}_2 , \tilde{\nu}_2 - \tilde{\nu}_1}(t_2)\right)\\
  & \times \frac{\prod_{j=1}^3 (1+ |it_1 + \alpha_j|)^{\frac{c_1-1}{2}} \prod_{j=1}^3 (1+|i t_2 - \alpha_j|)^{\frac{c_2-1}{2} }}{ (1+|t_1+t_2|)^{\frac{c_1+c_2-1}{2}} }y_1^{-c_1} y_2^{-c_2} dt_1\, dt_2.
  \end{split}
\end{displaymath}
We consider only the first summand in the first line of the preceding display. The second summand is similar. By a shift of variables, the first term equals
\begin{equation}\label{equals}
\begin{split}
     y_1^{1-c_1}y_2^{1-c_2} \int_{-\infty}^{\infty}\int_{-\infty}^{\infty} & w_{ 0 ,  3\tilde{\nu}_2}(t_1) w_{  0, 3\tilde{\nu}_1}(t_2)\frac{  (1+ |t_1 -3\tilde{\nu}_2|)^{\frac{c_1-1}{2}}  (1+ |t_1 |)^{\frac{c_1-1}{2}} (1+ |t_1 +3\tilde{\nu}_1|)^{\frac{c_1-1}{2}}  }{ (1+|t_1+t_2|)^{\frac{c_1+c_2-1}{2}} } \\
  &  \times  (1+| t_2  +3\tilde{\nu}_2|)^{\frac{c_2-1}{2} }(1+| t_2 |)^{\frac{c_2-1}{2} }(1+|t_2 - 3\tilde{\nu}_1|)^{\frac{c_2-1}{2} }   dt_1\, dt_2.
  \end{split}
\end{equation}
It is straightforward to estimate this expression. For convenience we provide the details. 
We recall our assumption \eqref{range} and split the $t_1, t_2$ integration into several ranges. Let $\tilde{\nu}_1 \leq R \leq \tilde{\nu_2}$, and define
\begin{displaymath}
\begin{split}
&\mathcal{I}_{-} := \{ t_1 \leq \tilde{\nu}_1\}, \quad  \mathcal{I}_{R} := \{R \leq t_1 \leq 2R\}, \quad \mathcal{I}_+ := \{ t_1 \geq 2\tilde{\nu}_2\},\\
& \mathcal{J}_{-} := \{t_2 \leq \tilde{\nu}_1\},  \quad \mathcal{J}_+ := \{t_2 \geq \tilde{\nu}_1\}. 
\end{split}
\end{displaymath}
We estimate the double integral in all 6 ranges for $c_1, c_2 > 0$:
\begin{displaymath}
\begin{split}
  \int_{\mathcal{I}_-} \int_{\mathcal{J}_-} \ll & \int_{-\infty}^{\tilde{\nu}_1} \int_{-\infty}^{\tilde{\nu}_1} \min(1, e^{\frac{\pi}{4} t_1}) \min(1, e^{\frac{\pi}{4} t_2})  ((1+\tilde{\nu}_1)(1+\tilde{\nu}_2))^{\frac{c_1+c_2}{2}-1}  \frac{(1+|t_1|)^{\frac{c_1-1}{2}}(1+|t_2|)^{\frac{c_2-1}{2}}}{(1+|t_1+t_2|)^{\frac{c_1+c_2-1}{2}}}     dt_2\, dt_1 \hspace{3cm} \\
&  \ll ((1+\tilde{\nu}_1)(1+\tilde{\nu}_2))^{\frac{c_1+c_2}{2}-1} (1+\tilde{\nu}_1)^{\frac{3}{2}};
  \end{split}
\end{displaymath}

 \begin{displaymath}
\begin{split}
 \int_{\mathcal{I}_-} \int_{\mathcal{J}_+}&  \ll  \int_{-\infty}^{\tilde{\nu}_1} \int_{\tilde{\nu}_1}^{\infty} \min(1, e^{\frac{\pi}{4} t_1})\min(1, e^{\frac{\pi}{4}(3\tilde{\nu}_1 - t_2)}) (1+\tilde{\nu}_1)^{ - \frac{1}
 {2}} (1 + \tilde{\nu}_2)^{\frac{c_1+c_2}{2}-1}  (1+ |t_1|)^{\frac{c_1-1}{2}}\hspace{3cm}  \\
 & \times (1+ |t_2 - 3\tilde{\nu}_1|)^{\frac{c_2-1}{2}} dt_2\, dt_1  \ll ((1+\tilde{\nu}_1)(1+\tilde{\nu}_2))^{\frac{c_1+c_2}{2}-1} (1+\tilde{\nu}_1)^{\frac{3}{2}};
 \end{split}
 \end{displaymath}

\begin{displaymath}
\begin{split}
  \int_{\mathcal{I}_R} \int_{\mathcal{J}_-} &\ll  \int_{R}^{2R} \int_{-\infty}^{\tilde{\nu}_1}  \min(1, e^{\frac{\pi}{4} t_2})  (1+\tilde{\nu}_1)^{\frac{c_2-1}{2}} (1+\tilde{\nu}_2)^{\frac{c_1+c_2}{2}-1} (1+R)^{\frac{c_1-c_2-1}{2}}  (1+|t_2|)^{\frac{c_2-1}{2}}       dt_2\, dt_1\hspace{3cm} \\
  & \ll (1+\tilde{\nu}_1) ^{c_2} (1+\tilde{\nu}_2)^{\frac{c_1+ c_2}{2}-1}  (1+R)^{\frac{c_1-c_2+1}{2}};  \end{split}
\end{displaymath}

 \begin{displaymath}
\begin{split}
 \int_{\mathcal{I}_R} \int_{\mathcal{J}_+}&  \ll  \int_{-R}^{2R} \int_{\tilde{\nu}_1}^{\infty} \min(1, e^{\frac{\pi}{4}(3\tilde{\nu}_1 - t_2)})   (1+\tilde{\nu}_1)^{\frac{c_2-1}{2}} (1+\tilde{\nu}_2)^{\frac{c_1+c_2}{2}-1} (1+R)^{\frac{c_1-c_2-1}{2}}    (1+ |t_2 - 3\tilde{\nu}_1|)^{\frac{c_2-1}{2}}     dt_2\, dt_1\hspace{3cm}\\
 &  \ll (1+\tilde{\nu}_1) ^{c_2} (1+\tilde{\nu}_2)^{\frac{c_1+ c_2}{2}-1}  (1+R)^{\frac{c_1-c_2+1}{2}};
 \end{split}
 \end{displaymath}

\begin{displaymath}
\begin{split}
  \int_{\mathcal{I}_+} \int_{\mathcal{J}_-} &\ll  \int_{2\tilde{\nu}_2}^{\infty} \int_{-\infty}^{\tilde{\nu}_1} \min(1, e^{\frac{\pi}{4}(3\tilde{\nu}_2 - t_1)})  \min(1, e^{\frac{\pi}{4} t_2})  (1+\tilde{\nu}_1)^{\frac{c_2-1}{2}} (1+\tilde{\nu}_2)^{\frac{c_1}{2}-1}  (1 + |t_1 - 3\tilde{\nu}_2|)^{ \frac{c_1-1}{2}} \hspace{3cm} \\
    &  \times  (1+|t_2|)^{\frac{c_2-1}{2}}   dt_2\, dt_1 \ll (1+\tilde{\nu}_1) ^{c_2} (1+\tilde{\nu}_2)^{c_1-\frac{1}{2}} ;
  \end{split}
\end{displaymath}

 \begin{displaymath}
\begin{split}
 \int_{\mathcal{I}_+} \int_{\mathcal{J}_+}&  \ll  \int_{2\tilde{\nu}_2}^{\infty} \int_{\tilde{\nu}_1}^{\infty}  \min(1, e^{\frac{\pi}{4}(3\tilde{\nu}_2 - t_1)})  \min(1, e^{\frac{\pi}{4} (3\tilde{\nu}_1-t_2)})  (1+\tilde{\nu}_1)^{\frac{c_2-1}{2}} (1+\tilde{\nu}_2)^{\frac{c_1}{2}-1}    (1 + |t_1 - 3\tilde{\nu}_2|)^{ \frac{c_1-1}{2}}   \hspace{3cm} \\
 &\times (1+ |t_2 - 3\tilde{\nu}_1|)^{\frac{c_2-1}{2}}  dt_2\, dt_1   \ll (1+\tilde{\nu}_1) ^{c_2} (1+\tilde{\nu}_2)^{c_1-\frac{1}{2}}.
  \end{split}
 \end{displaymath}
Combining all 6 previous bounds, and summing over dyadic numbers $R$, we obtain the bound of the proposition if $\nu_1, \nu_2 \in i \Bbb{R}$.

It remains to consider the situation \eqref{nu1}. The exponential factor does not change, but the fraction in the first line of \eqref{stir} becomes  now
\begin{equation}\label{except}
\begin{split}
 &  \frac{  (1+ |t_1 + \gamma|)^{\frac{c_1+\rho-1}{2}} (1+ |t_1 + \gamma|)^{\frac{c_1-\rho- 1}{2}} (1+ |t_1 - 2\gamma|)^{\frac{c_1-1}{2}}  }  { (1+|t_1+t_2|)^{\frac{c_1+c_2-1}{2}} } \\
&  \times (1+ |t_2 - \gamma|)^{\frac{c_2+\rho-1}{2}} (1+ |t_2 - \gamma|)^{\frac{c_2-\rho- 1}{2}} (1+ |t_2 + 2\gamma|)^{\frac{c_2-1}{2}}.  
\end{split}
\end{equation}
This is independent of $\rho$, and the same calculation goes through. \qed \\

We recall an important formula of Stade \cite{St} (cf.\ also \eqref{bar} and observe that Stade's definition \cite[(1.1)]{St} of the Whittaker function has $\nu_1$ and $\nu_2$ interchanged, and his Whittaker function is, up to Gamma-factors, twice our Whittaker function). 
\begin{prop}\label{prop2} 
Let $\nu_1, \nu_2, \mu_1, \mu_2, s \in \Bbb{C}$. We use the notation \eqref{defnu1} and \eqref{defnu2}, and define similarly $\mu_0$ and $\beta_1, \beta_2, \beta_3$ in terms of $\mu_1, \mu_2$. Then we have an equality of meromorphic functions in $s$:
\begin{displaymath}
\begin{split}
 & \int_0^{\infty} \int_0^{\infty} \overline{\tilde{W}_{\nu_1, \nu_2}(y_1, y_2)} \tilde{W}_{\mu_1, \mu_2}(y_1,y_2) (y_1^2y_2)^s \frac{dy_1\, dy_2}{(y_1y_2)^3} =  \frac{\pi^{3 -3s}  \prod_{j, k=1}^3\Gamma(\frac{s+ \alpha_j + \overline{\beta}_k}{2}) }{4 \Gamma(\frac{3}{2}s)\prod_{j=0}^2 |\Gamma(\frac{1}{2} + \frac{3}{2}i\Im \nu_j)| \prod_{j=0}^2 |\Gamma(\frac{1}{2} + \frac{3}{2}i \Im \overline{\mu}_j)|   }.
  \end{split}
\end{displaymath}
 \end{prop}
 
\section{Integrals over Whittaker functions} 
 
For our purposes it is convenient to consider the double Mellin transform of the product of $\tilde{W}_{\nu_1, \nu_2}$  with some rapidly decaying function. We are not aware of any explicit formula in the literature, but the next proposition gives an asymptotic result which is sufficient for our purposes. This is one of the key ingredients in this paper, and therefore we present all details of the lengthy proof. 

\begin{prop}\label{whitint} Let $A,  X_1, X_2 \geq 1$, $\tau_1, \tau_2 \geq 0$, and assume that $\tau_1+\tau_2$ is sufficiently large in terms of $A$.  Let $\nu_1, \nu_2\in \Bbb{C}$ satisfy \eqref{nu0} - \eqref{nu1} and in addition  $\Im \nu_1, \Im \nu_2 \geq 0$. Let   \begin{displaymath}
  t_1 = \tau_1 + 2\tau_2, \quad t_2 = 2\tau_1 + \tau_2. 
\end{displaymath} 
Fix two non-zero smooth functions $f_1, f_2 : (0, \infty) \rightarrow [0, 1]$ with compact support.    Let $\varepsilon > 0 $ and let
\begin{displaymath}
  I = I(\nu_1, \nu_2, t_1,t_2, X_1, X_2) = \left|\int_0^{\infty}\int_0^{\infty} \tilde{W}_{\nu_1, \nu_2}(y_1, y_2) f_1(X_1y_1) f_2(X_2y_2)  y_1^{it_1}y_2^{it_2} \frac{dy_1\, dy_2}{(y_1y_2)^3}\right|. 
\end{displaymath}
 If $X_1 = X_2 = 1$, then
\begin{equation}\label{bound}
\begin{split}
 I   \ll   (|\Im \nu_1 - \tau_1| + |\Im\nu_2 - \tau_2|)^{-A} \prod_{j=0}^2 (1+|\nu_j|)^{-\frac{1}{2}}      \end{split}
\end{equation}
Moreover, there is a constant $c$ depending on $f_1, f_2$ such that the following holds: if 
\begin{equation}\label{fit}
  X_1 = X_2 = 1, \quad \tau_1, \tau_2 \gg   1, \quad   |\Im \nu_1 - \tau_1| \leq c, \quad |\Im \nu_2 - \tau_2| \leq c,
\end{equation}
then 
\begin{equation}\label{asymp}
  I  \asymp   \prod_{j=0}^2 (1+|\nu_j|)^{-1/2};
\end{equation}
and if $\nu_1, \nu_2$ are given by \eqref{nontemp} and in addition
\begin{equation}\label{nontempasump}
 X_2 \gg 1, \quad |\rho| \geq \varepsilon, \quad  \gamma \gg  (X_1X_2)^{\varepsilon}, \quad |\gamma - \tau_2| \leq c, \quad |\tau_1| \leq c, 
\end{equation}
then
 \begin{equation}\label{asymp1}
  I  \asymp X_1 X_2^{1+ |\rho|} \prod_{j=0}^2 (1+|\nu_j|)^{-1/2}. 
\end{equation}
All implied constant depend at most on $A, \varepsilon, f_1, f_2$ and the sign $\gg$ should be interpreted as ``up to a sufficiently large constant".  
 \end{prop}

\textbf{Remark 2:} This should be roughly  interpreted as follows: given $t_1, t_2, X_1, X_2$ as above, $I$ as a function of $\nu_1, \nu_2 \in \{z \in \Bbb{C} : |\Im z| \leq 1/2\}$ is under some technical assumptions a function with a bump at $\Im \nu_1 = \tau_1$ and $\Im \nu_2 = \tau_2$ of size $X_1  X_2^{1+\max_j|\Re \alpha_j|}(\nu_0\nu_1\nu_2)^{-1/2}$ with rapid decay away from this point.  Most of the time we will put $X_1=X_2 = 1$. Only if we need a test function that blows up at exceptional eigenvalues we will choose $X_2$ to be large. The asymmetry  in $X_1, X_2$ in \eqref{nontempasump} and \eqref{asymp1} is due to the special choice \eqref{nontemp}.  \\

\textbf{Proof.}  By Parseval's formula and \eqref{doubleMellin} the double integral in question equals
\begin{equation}\label{parse}
\begin{split}
 & \int_{(1/2)} \int_{(1/2)}  
  \frac{\widehat{f_1}(-1+it_1-u_1) \widehat{f_2}(-1+it_2-u_2)
 \prod_{j=1}^3 \Gamma(\frac{1}{2}(u_1 + \alpha_j)) \prod_{j=1}^3 \Gamma(\frac{1}{2}(u_2 - \alpha_j)) }{ 4\pi^{u_1+u_2 - \frac{3}{2}} \Gamma(\frac{1}{2}(u_1+u_2))  \prod_{j=0}^2  |\Gamma\left(\frac{1}{2} + \frac{3}{2}i \Im \nu_j\right)|X_1^{-1+it_1-u_1}X_2^{-1+it_2-u_2 }} \frac{du_1\, du_2}{(2\pi i)^2}.
  \end{split}
\end{equation}
Let us first assume that
\begin{equation}\label{assump}
  |\nu_1| + |\nu_2| \leq \frac{1}{100}(\tau_1 + \tau_2). 
\end{equation}
In this case the conditions \eqref{fit} and \eqref{nontempasump} are void, so we only need to show \eqref{bound} and take $X_1 = X_2 = 1$. We  
apply Stirling's formula to the $\Gamma$-quotient. We argue as in the proof of Proposition \ref{prop1}, see \eqref{exp} and the surrounding discussion. The exponential part is given by
\begin{displaymath}
  \exp\left(\frac{3\pi}{4}\sum_{j=0}^2|\Im \nu_j|   - \frac{\pi}{4}\sum_{j=1}^3|\Im (u_1 +  \alpha_j)| - \frac{\pi}{4}\sum_{j=1}^3|\Im (u_2 -  \alpha_j)| + \frac{\pi}{4} |\Im (u_1+u_2)| \right). 
\end{displaymath}
As before let us write $\tilde{\nu}_j = \Im \nu_j$ and assume without loss of generality \eqref{range}. Using  \eqref{case1} and \eqref{case2} together with the rapid decay of $\widehat{f_1}$ and $\widehat{f_2}$ it is easy to see that by our present assumption  \eqref{assump} we can bound $I$ by 
\begin{displaymath}
 \ll_{A, f_1, f_2}  (t_1+t_2)^{-A}. 
\end{displaymath}
In the range \eqref{assump} this is acceptable for \eqref{bound}.  \\

Let us now assume 
\begin{equation}\label{newassump}
  |\nu_1| + |\nu_2| \geq \frac{1}{100}(\tau_1 + \tau_2). 
\end{equation} 
We want to shift the two contours in \eqref{parse} to $-\infty$. To check convergence, we first shift the $u_1$-integral to $\Re u_1 = -2A-1$ for some large integer $A$. We observe that 
\begin{displaymath}
  \widehat{f}(s) \ll_{B, f} |s|^{-B} C^{\Re s+B}
\end{displaymath}
for $\Re s > 0$, any $B \geq 0$ and some constant $C > 0$ depending only on $f$ (one can take $C := \sup\{x > 0 \mid f(x) \not= 0\}$). We also recall that the reflection formula for the Gamma-function implies the uniform bound
\begin{displaymath}
  \Gamma(-s) \ll e^{-\frac{\pi}{2}|\Im s|} \left(\frac{|s|}{e}\right)^{-\Re s - \frac{1}{2}} 
\end{displaymath}
for $\Re s > 0$, $\min_{n \in \Bbb{Z}}(\Re s - n) > 1/50$.  It is now easy to see that the remaining $u_1$-integral for $A\rightarrow \infty$ vanishes, and we are left with the sum over the residues. Next we shift in the same way the $u_2$-integral to $-\infty$, and express \eqref{parse} as an absolutely convergent double sum over residues. 
 Let us first assume that  $\nu_1 \not=0$ so that $\alpha_1, \alpha_2, \alpha_3$ are pairwise distinct. For $j \in \{1, 2, 3\}$ we denote by $k, l$ two integers such that $\{j, k, l\} = \{1, 2, 3\}$. Similarly for $r \in \{1, 2, 3\}$ let $s, t$ be such that $\{r, s, t\} = \{1, 2, 3\}$.  Then \eqref{parse} equals  \begin{equation}\label{mainterm}
\begin{split}
 &  \underset{\substack{1 \leq j, r \leq 3\\ j \not = r}}{\sum\sum}   
 \frac{\Gamma(\frac{-\alpha_j+\alpha_k}{2})\Gamma(\frac{-\alpha_j+ \alpha_l}{2}) \Gamma(\frac{\alpha_r - \alpha_s}{2})\Gamma(\frac{\alpha_r - \alpha_t}{2}) }{\prod_{j=0}^2 |\Gamma(\frac{1}{2} + \frac{3}{2}  i \Im \nu_j)|\, 4\pi^{-\alpha_j+\alpha_r-\frac{3}{2}} \Gamma(\frac{-\alpha_j + \alpha_r}{2}) } \\
& \times  \sum_{n=0}^{\infty} \sum_{m=0}^{\infty} \frac{4(-1)^{n+m}(\frac{-\alpha_j + \alpha_r}{2}-1)_{n+m}  }{n!m! (\frac{-\alpha_j+\alpha_k}{2}-1)_n(\frac{-\alpha_j+ \alpha_l}{2}-1)_n (\frac{\alpha_r - \alpha_s}{2}-1)_m (\frac{\alpha_r - \alpha_t}{2}-1)_m }\\
& \quad\quad \times \frac{ \widehat{f_1}(-1+\alpha_j+it_1 + 2n)\widehat{f_2}(-1-\alpha_r + it_2 +2m)}{ X_1^{-1+\alpha_j+it_1+2n}X_2^{-1-\alpha_r+it_2 + 2m}}.
\end{split}
\end{equation}
  
We can   bound the second and third  line of \eqref{mainterm}  by
\begin{equation}\label{crude}
\sum_{n, m}   \ll_{A, f_1, f_2}  X_1^{1-\Re \alpha_j} X_2^{1+\Re\alpha_r} (1+|\alpha_j + it_1|)^{-A} (1+|\alpha_r-it_2|)^{-A}.
\end{equation}
Here we have used  that by \eqref{defnu2} we  have the following equality of multisets:
\begin{equation}\label{multiset}
  \{-\alpha_j + \alpha_k, -\alpha_j +\alpha_l, \alpha_r - \alpha_s, \alpha_r - \alpha_t\} \setminus \{-\alpha_j + \alpha_r\} = \{\pm 3\nu_0, \pm 3\nu_1, \pm 3 \nu_2\}
\end{equation}
for a certain choice of signs (depending on $j, r$) whenever $j \not= r$.  In addition we see that   we have in the special case $j=3$, $r\in \{1, 2\}$ (that is, $\alpha_j = -\nu_1-2\nu_2$, $\alpha_r = 2\nu_1 + \nu_2$ or $-\nu_1 + \nu_2$)
\begin{equation}\label{without0}
  \sum_{(n, m) \not= (0, 0)} \ll_{A,   f_1, f_2}  X_1^{1-\Re \alpha_j} X_2^{1+\Re\alpha_r} \frac{(1+|\alpha_j + it_1|)^{-A} (1+|\alpha_r-it_2|)^{-A}.
}{\min\left( (1+ |\nu_1|)X_2^2, (1+ |\nu_2|)X_1^2)\right)}.
\end{equation}
 We will now carefully analyze all 6 terms $1 \leq \alpha_j, \alpha_r \leq 3$, $ j \not= r$ in the main term under the assumption $\tau_1, \tau_2 \geq 0,  \Im \nu_2 \geq  \Im \nu_1 \geq 0$ and show that they all satisfy the bound \eqref{bound}. Moreover, under the assumption \eqref{fit}, the term $j=3$, $r=1$ is of order of magnitude \eqref{asymp} and dominates all other terms. Similarly we will show \eqref{asymp1}.  We will first make the extra assumption
\begin{displaymath}
  |\nu_1| \geq \varepsilon. 
\end{displaymath}
This ensures that $\alpha_1, \alpha_2, \alpha_3$ are not too close together (note that $|\nu_2|$ must be large  by \eqref{newassump}).  
By Stirling's formula, \eqref{multiset} and \eqref{nontemp} (in the non-tempered case), 
\begin{equation}\label{stir1}
\begin{split}
 & \left|\frac{\Gamma(\frac{-\alpha_j+\alpha_k}{2})\Gamma(\frac{-\alpha_j+ \alpha_l}{2}) \Gamma(\frac{\alpha_r - \alpha_s}{2})\Gamma(\frac{\alpha_r - \alpha_t}{2}) }{\prod_{j=0}^2 |\Gamma(\frac{1}{2} + \frac{3}{2}  i \Im \nu_j)| \Gamma(\frac{-\alpha_j + \alpha_r}{2})} \right| \begin{cases}
 \ll_{\varepsilon}  \prod_{n=0}^2 (1+|\nu_n|)^{O(1)}, & \{j, r\} = \{1, 2\} \text{ and }   \nu_1 \in \Bbb{R}, \\
 \asymp_{\varepsilon} \prod_{n=0}^2 (1+|\nu_n|)^{ -\frac{1}{2}}, & \text{otherwise}.  
  \end{cases}
 \end{split}
\end{equation}
 The non-negativity of $f_1, f_2$ implies that the absolute values of the Mellin transforms $|\widehat{f_1}(-1+\alpha_3 + it_1)|, |\widehat{f_2}(-1-\alpha_1+ i t_2)|$ are bounded from below  in the range \eqref{fit} and \eqref{nontempasump} if $c$ is sufficiently small. Combining \eqref{stir1} and \eqref{without0}, we see that under the assumption \eqref{fit} the term $j=3$, $r=1$ satisfies \eqref{asymp}. Note that \eqref{fit} forces $|\nu_1|, |\nu_2|$ to be sufficiently large and excludes \eqref{nu1}.    Similarly, under the assumption \eqref{nontempasump} the term $j=3$, $r=1$ (if $\rho > 0$) or the term $j=3$, $r=2$ (if $\rho < 0$) satisfies \eqref{asymp1}, whereas the other term is of smaller order of magnitude. This is also consistent with \eqref{bound}. \\

It remains to show that all other terms satisfy \eqref{bound}, and are of lesser order of magnitude than \eqref{asymp} and \eqref{asymp1} under the respective conditions. 
 Under the assumption \eqref{fit} all   5 terms $(j, r) \not= (3, 1)$ satisfy $| \alpha_j + it_1| + |\alpha_r - it_2| \geq 3 \min(|\Im \nu_1|, |\Im \nu_2|) + O(1)$  and can therefore be bounded by (recall \eqref{range})
\begin{displaymath}
  \ll_{\varepsilon,  A,   f_1, f_2}  |\nu_1|^{-A}   \prod_{n=0}^2 (1+|\nu_n|)^{-\frac{1}{2}}
\end{displaymath}
which is dominated by \eqref{asymp}. Similarly, under the assumption \eqref{nontempasump} the 4 terms $(\alpha_j,  \alpha_r) \not\in\{ (3, 1), (3, 2)\}$ satisfy $| \alpha_j + it_1| + |\alpha_r - it_2| \geq 3 |\Im \nu_2| + O(1)$  
and can therefore be bounded by
\begin{displaymath}
  \ll_{\varepsilon,  A,   f_1, f_2}  |\nu_2|^{-A}   \prod_{n=0}^2 (1+|\nu_n|)^{-\frac{1}{2}}
\end{displaymath}
which is again dominated by \eqref{asymp1}.  We proceed now to show \eqref{bound} for $X_1=X_2 = 1$.  It follows from \eqref{crude} and \eqref{stir1} that all 6 terms   $\alpha_j  \not= \alpha_r$ contribute
\begin{displaymath}
  \ll_{\varepsilon, A, f_1, f_2}    ((1+ | \alpha_j +i t_1|)( 1+ | \alpha_r - it_2|))^{-A}   \begin{cases}
   \prod_{n=0}^2 (1+|\nu_n|)^{O(1)}, & \{j, r\} = \{1, 2\} \text{ and }   \nu_1 \in \Bbb{R}, \\
   \prod_{n=0}^2 (1+|\nu_n|)^{ -\frac{1}{2}}, & \text{otherwise}.  
  \end{cases}
\end{displaymath}
to the main term. This is in agreement with \eqref{bound} if we can show
\begin{displaymath}
  |\Im \alpha_j + t_1| + |\Im \alpha_r - t_2 | \geq  \frac{1}{2}\left(|\tilde{ \nu}_1 - \tau_1| + |\tilde{\nu}_2 - \tau_2|\right).
\end{displaymath}
This is is clear for $j=1$ \emph{or} $r=3$ by the positivity assumption $\tilde{ \nu}_1, \tilde{\nu}_2, \tau_1, \tau_2 \geq 0$ (recall the notation $\tilde{\nu}_j = \Im \nu_j$), even without the factor $1/2$.  We check the other 3 cases. In the case $j=3, r=1$ we have again the stronger inequality
\begin{displaymath}
  |-\tilde{\nu}_1 - 2\tilde{\nu}_2 + \tau_1 + 2 \tau_2| + |2\tilde{ \nu}_1 + \tilde{\nu}_2 - 2\tau_1- \tau_2| \geq |\tau_1 -\tilde{\nu}_1| + |\tau_2 - \tilde{\nu}_2|
\end{displaymath}
which follows from the easy to check inequality $|a| + |b| \leq |a+2b| + |b+2a|$. In the case $j=3$, $r=2$ we need to show
\begin{displaymath}
  |-\tilde{ \nu}_1 -2\tilde{ \nu}_2 + \tau_1 + 2 \tau_2| + |-\tilde{ \nu}_1 +\tilde{\nu}_2 - 2\tau_1- \tau_2| \geq \frac{1}{2}\left(|\tau_1 -\tilde{\nu}_1| + |\tau_2 - \tilde{ \nu}_2|\right).
\end{displaymath}
If $\tau_1 -\tilde{ \nu}_1$ and $\tau_2-\tilde{ \nu}_2$ are of the same sign, the first term dominates the right hand side; if $\tau_1 - \tilde{\nu}_1 \leq 0, \tau_2 - \tilde{ \nu}_2 \geq 0$, the second term dominates the right hand side; in either case we do not need the factor $1/2$. Finally if  $\tau_1 - \tilde{\nu}_1 \geq 0, \tau_2 - \tilde{\nu}_2 \leq 0$, we distinguish the two cases $\tau_1 - \tilde{\nu}_1$ greater or smaller than  $\tilde{\nu}_2 - \tau_2$: in the former case the second term dominates the right hand side, because  $|-\tilde{\nu}_1 +\tilde{\nu}_2 - 2\tau_1- \tau_2|  \geq 2(\tau_1 - \tilde{\nu}_1) - (\tilde{\nu}_2 - \tau_2)$, and in the latter the first term dominates the right hand side, because $|-\tilde{ \nu}_1 -2\tilde{ \nu}_2 + \tau_1 + 2 \tau_2| \geq 2(\tilde{\nu}_2 -\tau_2) - (\tau_1 - \tilde{\nu}_1)$. Finally the case $j=2, r=1$ amounts to showing 
\begin{displaymath}
  |-\tilde{ \nu}_1 +\tilde{\nu}_2 + \tau_1 + 2 \tau_2| + |2\tilde{ \nu}_1 +\tilde{\nu}_2 - 2\tau_1- \tau_2| \geq \frac{1}{2}\left(|\tau_1 -\tilde{\nu}_1| + |\tau_2 - \tilde{\nu}_2|\right)
\end{displaymath}
which can be seen as above after interchanging  indices.\\ 

Finally we need to treat the case $0 \not= |\nu_1| < \varepsilon$ and $|\nu_2| \gg |\tau_1| + |\tau_2|$. Here the condition \eqref{nontempasump} is empty, and if  $\tau_1, \tau_2$ are sufficiently large, the condition \eqref{fit} is also empty, so we only need to show the upper bound \eqref{bound} for $X_1 = X_2 = 1$. 
We return to \eqref{mainterm} and partition the 6 terms $(j, r)$ into three pairs
\begin{displaymath}
   \{(3, 2), (3, 1)\}, \quad \{(2, 3), (1, 3)\}, \quad \{(2, 1), (1, 2)\}. 
\end{displaymath}
The contribution of the first pair is
\begin{displaymath}
\begin{split}
& \sum_{n, m = 0}^{\infty}  \frac{(-1)^{n+m}}{n! m!  } \frac{\Gamma(\frac{3\nu_2}{2} - n)\Gamma(\frac{3\nu_0}{2} - n)\Gamma(\frac{3\nu_2}{2} - m)\widehat{f_1}(-1-\nu_1-2\nu_2 + s_1 + 2n)}{\Gamma(\frac{3\nu_2}{2} - n - m) \prod_{j=0}^2 |\Gamma(\frac{1}{2} + \frac{3}{2}  i \Im \nu_j)| 4\pi^{-\frac{3}{2} +\nu_1+3\nu_2} }\\
 & \times   \left(\frac{\Gamma(\frac{-3\nu_1}{2} - m)\widehat{f_2}(-1+\nu_1 - \nu_2+it_2+2m)}{\pi^{-\nu_1} }  + \frac{\Gamma(\frac{3\nu_1}{2} - m)\widehat{f_2}(-1-2\nu_1 - \nu_2+it_2+2m)}{\pi^{2\nu_1} }\right). 
\end{split}
\end{displaymath}
For $|\nu_1| < \varepsilon$ the second line can be bounded by the mean value theorem. 
 Then we use the functional equation $s \Gamma(s) = \Gamma(s+1)$ of the Gamma-function in connection with Stirling's formula as before and bound the preceding display by
\begin{displaymath}
\begin{split}
& \ll_{A, f_1, f_2}   (| \alpha_1 +i t_1| +   | \nu_2 +i t_2|)^{-A}  \prod_{n=0}^2 (1+|\nu_n|)^{-\frac{1}{2}}
 \end{split}
\end{displaymath} 
and argue as before. The same argument with different indices works for the pair $\{(2, 3), (1, 3)\}$. The  last pair is only a small variation; its contribution is given by 
\begin{displaymath}
\begin{split}
&\sum_{n, m =0}^{\infty} \frac{(-1)^{n+m}\pi^{3/2}}{n!m!  4   \prod_{j=0}^2|\Gamma(\frac{1}{2} + \frac{3}{2}  i \Im \nu_j)| }\\
\times &  \left(\frac{\Gamma(\frac{-3\nu_2}{2} - n)\Gamma(\frac{3\nu_1}{2} - n)\Gamma(\frac{3\nu_0}{2} - m) \Gamma(\frac{3\nu_1}{2}-m)\widehat{f_1}(-1-\nu_1+\nu_2 + it_1 + 2n)\widehat{f_2}(-1- 2\nu_1 - \nu_2+it_2 + 2m)}{\pi^{3\nu_1} \Gamma(\frac{3\nu_1}{2} - n - m)  } \right.\\
 &  +\frac{\Gamma(\frac{-3(\nu_2+\nu_1)}{2} - n) \Gamma(\frac{-3\nu_1}{2} - n)\Gamma(\frac{3(\nu_0-\nu_1)}{2} - m)\Gamma(\frac{-3\nu_1}{2}-m)}{\pi^{-3\nu_1} \Gamma(\frac{-3\nu_1}{2} - n - m) }  \\
& \left.  \quad\quad\quad\quad\quad\quad \times   \widehat{f_1}(-1+2\nu_1+\nu_2 + it_1 + 2n)\widehat{f_2}(-1 + \nu_1 - \nu_2+it_2 + 2m)    \right). 
  \end{split}
\end{displaymath}
For small $\nu_1$,  this can again be estimated by the mean value theorem giving the crude bound
\begin{displaymath}
  \ll_{A, f_1, f_2} (| \nu_2 + it_1| + | \nu_2 -i t_2|)^{-A}  (1+|\nu_2|)^{O(1)}  
\end{displaymath}
which is admissible for \eqref{bound}. 
This completes the  proof of the proposition under the   additional assumption that $\alpha_1, \alpha_2, \alpha_3$ are pairwise distinct, that is $\nu_1 \not= 0$.  The case $\nu_1 = 0$ follows by continuity. \qed\\

An inspection of the proof, in particular \eqref{mainterm} -- \eqref{stir1}, shows that for $\tau_1$, $\tau_2$ sufficiently large one has
\begin{equation}\label{asympwhit}
 | I(\nu_1, \nu_2, t_1, t_2, 1, 1)|^2 \sim \frac{(2\pi)^3}{3^3 |\nu_0\nu_1\nu_2|} |\widehat{f_1}(-1+i\tau_1 - \nu_1 + 2i\tau_2 - 2\nu_2)\widehat{f_2}(-1+ 2 i\tau_1 -2 \nu_1 + i\tau_2 - \nu_2)|^2
\end{equation}
for $\nu_1, \nu_2 \in i\Bbb{R}$ in a neighbourhood of $i\tau_1$, $i\tau_2$, respectively, and it is very small outside this region.

\section{Maass forms}

 Let $\Gamma = SL_3(\Bbb{Z})$. We denote by 
\begin{displaymath}
  P  =  \left(\begin{matrix} \ast & \ast & \ast\\ \ast& \ast & \ast \\ 0 & 0 & \ast\end{matrix} \right) \subseteq \Gamma, \quad  P_1  =  \left(\begin{matrix} \ast & \ast & \ast\\ \ast& \ast & \ast \\ 0 & 0 & 1\end{matrix} \right) \subseteq \Gamma
\end{displaymath}  
the maximal parabolic subgroup, and by
\begin{displaymath}
 U_3 =  \left(\begin{matrix} 1 & \ast & \ast\\ 0& 1 & \ast \\ 0 & 0 & 1\end{matrix} \right) \subseteq \Gamma
\end{displaymath}  
the standard unipotent group. Analogously, let $U_2 := \left(\begin{matrix} 1 & \ast  \\ 0& 1 \end{matrix} \right) \subseteq SL_2(\Bbb{Z})$.

A Maa{\ss} cusp form $\phi : \Gamma\backslash \mathfrak{h}^3 \rightarrow \Bbb{C}$ with spectral parameters $\nu_1, \nu_2$ (that is, of type $(1/3 + \nu_1, 1/3 +  \nu_2)$ in the notation of \cite{Go}) for the group $\Gamma$ has a Fourier expansion of the type
\begin{equation}\label{defphi}
  \phi(z) = \sum_{m_1 =1}^{\infty} \sum_{m_2 \not=0} \frac{A_{\phi}(m_1, m_2)}{|m_1m_2|} \sum_{\gamma  \in U_2\backslash SL_2(\Bbb{Z})}\mathcal{W}^{\text{sgn}(m_2)}_{\nu_1, \nu_2}\left(\left(\begin{smallmatrix} |m_1m_2| & &\\ & m_1 & \\  & & 1\end{smallmatrix}\right) \left(\begin{smallmatrix}\gamma  & \\ &  1\end{smallmatrix}\right) z\right) c_{\nu_1, \nu_2}
\end{equation}  
with $\mathcal{W}^{\pm}_{\nu_1, \nu_2}$ as in \eqref{generalWhit} and $c_{\nu_1, \nu_2}$ as in \eqref{const}. The Fourier coefficients are given by
\begin{displaymath}
  \int_0^1\int_0^1 \int_0^1 \phi(z) e(-m_1x_1-m_2x_2) dx_1 \, dx_2 \, dx_3 = \frac{A_{\phi}(m_1, m_2)}{|m_1m_2|} \tilde{W}_{\nu_1, \nu_2}(m_1y_1, |m_2|y_2). 
\end{displaymath}
We have
\begin{displaymath}
  A_{\phi}(m_1, m_2) = A_{\phi}(m_1, -m_2),
\end{displaymath}
see \cite[Proposition 6.3.5]{Go}. Hence one can alternatively write the Fourier expansion as a sum over $m_1, m_2 \geq 1$, $\gamma \in U_2 \backslash GL_2(\Bbb{Z})$. We will use this observation in the proof of Lemma \ref{lem1}. 

  It is expected that $\nu_1, \nu_2$ are imaginary, but we certainly know that \eqref{nu0} -- \eqref{nu1} hold.  If $\phi$ is an eigenfunction of the Hecke algebra (see \cite[Section 6.4]{Go}), we define its $L$-function by $L(\phi, s) := \sum_{m} A_{\phi}(1,m) m^{-s}$, and the Rankin-Selberg $L$-function by
\begin{equation}\label{RS}
  L(\phi \times \tilde{\phi}, s) :=  \zeta(3s) \sum_{m_1, m_2} \frac{|A_{\phi}(m_1, m_2)|^2}{m_1^{2s} m_2^{s}}.  
\end{equation}
It follows from   \cite[Theorem 2]{Li1} or \cite[Corollary 2]{Brm} that the coefficients are essentially bounded on average, uniformly in $\nu$:
\begin{equation}\label{ranksel}
  \sum_{m \leq x} |A_{\phi}(m, 1)|^2 \ll x (x(1+|\nu_1|+|\nu_2|))^{\varepsilon}.
\end{equation}
The space of cusp forms is equipped with an inner product
\begin{displaymath}
   \langle \phi_1, \phi_2 \rangle := \int_{\Gamma \backslash \mathfrak{h}^3}\phi_1(z) \overline{\phi_2(z)} dx_1\, dx_2\, dx_3\, \frac{dy_1\, dy_2}{(y_1y_2)^3}.  
\end{displaymath}
It is known that $L(\phi \times \tilde{\phi}, s)$ can be continued holomorphically to $\Bbb{C}$ with the exception of a simple pole at $s=1$ whose residue is proportional to $\|\phi \|^2$ \cite[Theorem 7.4.9]{Go}. The proportionality constant is given in the next lemma. 
\begin{lemma}\label{lem1} 
For a Hecke eigenform $\phi$ as in \eqref{defphi} with $A_{\phi}(1, 1) = 1$ we have
 $ \|\phi \|^2 \asymp   
  \underset{s=1}{\text{\rm res}}   L(\phi \times \tilde{\phi}, s)$. 
 \end{lemma}

\textbf{Remark 3:} This lemma shows that the normalization of the Whittaker functions $\tilde{W}_{\nu_1, \nu_2}$ is well chosen in the sense that an arithmetically normalized cusp form $A_{\phi}(1, 1) = 1$ should roughly have norm 1. The main point is that $\tilde{W}_{\nu_1, \nu_2}$ has roughly norm 1 with respect to the inner product  \begin{displaymath}
  (f, g) := \int_0^{\infty} \int_0^{\infty} f(y_1, y_2) \overline{g(y_1, y_2)} \det\left(\begin{smallmatrix} y_1y_2 & & \\ & y_1 & \\ & & 1\end{smallmatrix}\right) \frac{dy_1 \, dy_2}{(y_1y_2)^3}.
\end{displaymath}\\
 
\textbf{Proof.} This is standard Rankin-Selberg theory. We use the maximal parabolic Eisenstein series 
\begin{displaymath}
  E(z, s; \textbf{1}) := \sum_{\gamma \in P \backslash \Gamma} \det(\gamma z)^s = \frac{1}{2} \sum_{\gamma \in P_1 \backslash \Gamma} \det(\gamma z)^s  , \quad \Re s > 1.
\end{displaymath}  
   It follows from \eqref{res} below that
\begin{displaymath}
  \| \phi \|^2 =  \frac{3\zeta(3)}{2\pi}
\underset{s=1}{\text{res}}\langle \phi, \phi  E(., \bar{s}, \textbf{1}) \rangle. 
\end{displaymath}
We follow the unfolding argument of \cite[p.\ 227-229]{Go} and \cite[Section 3]{Fr}. 
Unfolding the Eisenstein series, we see 
\begin{displaymath}
  \langle \phi, \phi E(., \bar{s}, \textbf{1}) \rangle = \frac{1}{2} \int_{P_1\backslash \mathfrak{h}^3} |\phi(z)|^2 (y_1^2y_2)^s dx_1\, dx_2\, dx_3\, \frac{dy_1\, dy_2}{(y_1y_2)^3}.
\end{displaymath}
Let $\mathcal{F}$ denote a fundamental domain for $\left\{\left(\begin{smallmatrix}1 & & \ast\\ & 1 & \ast \\ & & 1\end{smallmatrix}\right)\right\} \backslash \mathfrak{h}^3$, and let $\widetilde{GL}_2(\Bbb{Z}) := \left\{\left(\begin{smallmatrix} \gamma & \\ & 1\end{smallmatrix}\right) \mid \gamma \in GL_2(\Bbb{Z})\right\} \subseteq GL_3(\Bbb{Z})$. Then $P_1\backslash \mathfrak{h}^3$ is in 2-to-1 correspondence with $\widetilde{GL}_2(\Bbb{Z}) \backslash \mathcal{F}$. Inserting the Fourier expansion of one factor and unfolding once again, we obtain
\begin{displaymath}
\begin{split}
  \langle \phi, \phi E(., \bar{s}, \textbf{1}) \rangle & = \sum_{m_1=1}^{\infty} \sum_{m_2 =1}^{\infty} \frac{|A_{\phi}(m_1, m_2)|^2}{|m_1m_2|^2} \int_{-\infty}^{\infty} \int_{-\infty}^{\infty}  |\tilde{W}_{\nu_1, \nu_2}(m_1y_1, |m_2|y_2)|^2 (y_1^2y_2)^{s} \frac{dy_1\, dy_2}{(y_1y_2)^3}\\
 & = \frac{L(\phi \times \tilde{\phi}, s)  }{\zeta(3s)}   \int_{-\infty}^{\infty} \int_{-\infty}^{\infty}  |\tilde{W}_{\nu_1, \nu_2}(y_1, y_2)|^2 (y_1^2y_2)^{s} \frac{dy_1\, dy_2}{(y_1y_2)^3}. 
 \end{split} 
\end{displaymath}
  The lemma follows now easily from Stade's formula. \qed \\ 
 
We are now ready to prove \eqref{resiweak}.

\begin{lemma}\label{lem2} For an arithmetically normalized Hecke-Maa{\ss}  cusp form $\phi$ with spectral parameters $\nu_1, \nu_2$ as above we have
\begin{displaymath} 
  \left((1+|\nu_0|) (1+|\nu_1|) (1+|\nu_2|)\right)^{-1} \ll \| \phi\|^2 \ll_{\varepsilon}   \left((1+|\nu_0|) (1+|\nu_1|) (1+|\nu_2|)\right)^{\varepsilon}
\end{displaymath}
for any $\varepsilon > 0$. 
\end{lemma} 
 
\textbf{Proof.} We conclude from Lemma \ref{lem1} that  as in \eqref{ranksel} the upper bound follows directly from \cite[Theorem 2]{Li1} or \cite[Corollary 2]{Brm}.    
We proceed to prove the lower bound. The idea is taken from \cite[Lemma 4]{DI2}. We can assume that one of $\nu_1, \nu_2$ is sufficiently large. Since $A_{\phi}(1, 1)= 1$, we have
\begin{displaymath}
  \tilde{W}_{\nu_1, \nu_2}(y_1, y_2) = \int_0^1\int_0^1\int_0^1 \phi(z) e(-x_1-x_2)  dx_1 \, dx_2 \, dx_3.
\end{displaymath}
for any $y_1, y_2 > 0$. By Cauchy-Schwarz we get
\begin{displaymath}
   | \tilde{W}_{\nu_1, \nu_2}(y_1, y_2)|^2 \leq  \left(\int_0^1\int_0^1\int_0^1 |\phi(z)|^2   dx_1 \, dx_2 \, dx_3\right)^{1/2} |\tilde{W}_{\nu_1, \nu_2}(y_1, y_2)|.
\end{displaymath}
  Integrating this inequality and using Cauchy-Schwarz again, we find
\begin{displaymath}
\begin{split}
& \int_1^{\infty} \int_1^{\infty}  | \tilde{W}_{\nu_1, \nu_2}(y_1, y_2)|^2 (y_1^2 y_2)^{1/2}  \frac{dy_1\, dy_2}{y_1^3y_2^3} \\
 &\leq  \left(\int_1^{\infty} \int_1^{\infty} \int_0^1\int_0^1\int_0^1 |\phi(z)|^2    \frac{dx_1 \, dx_2 \, dx_3 \, dy_1 \, dy_2}{(y_1y_2)^3}\right)^{1/2} \left(\int_0^{\infty} \int_0^{\infty} |\tilde{W}_{\nu_1, \nu_2}(y_1, y_2)|^2 y_1^2y_2 \frac{dy_1\, dy_2}{(y_1y_2)^3}\right)^{1/2}.
\end{split} 
\end{displaymath}
Since $[1, \infty)^2 \times [0, 1]^3$ is contained in a fundamental domain for $SL_3(\Bbb{Z}) \backslash \mathfrak{h}^3$ (see e.g.\ \cite{Gr}), we obtain together with Proposition \ref{prop2} that the right hand side is
\begin{displaymath}
 \ll \| \phi \|  \left(  \frac{ \prod_{j, k=1}^3 \Gamma(\frac{1+ \alpha_j +  \overline{\alpha}_k}{2}) }{ \prod_{j=0}^2 |\Gamma(\frac{1}{2} + \frac{3}{2}i\Im\nu_j)|^2 } \right)^{1/2} \asymp \|\phi\|. 
\end{displaymath}
 The left hand side is
\begin{displaymath}
  \geq  \int_0^{\infty} \int_0^{\infty}  | \tilde{W}_{\nu_1, \nu_2}(y_1, y_2)|^2 (y_1^2 y_2)^{1/2}  \frac{dy_1\, dy_2}{y_1^3y_2^3} -  \int_0^{\infty} \int_0^{\infty}  | \tilde{W}_{\nu_1, \nu_2}(y_1, y_2)|^2 (y_1^2 y_2)^{1/4}  \frac{dy_1\, dy_2}{y_1^3y_2^3}. 
\end{displaymath}
Again by Proposition \ref{prop2}, this is
\begin{displaymath}
\begin{split}
& \asymp   \left((1+|\nu_0|) (1+|\nu_1|) (1+|\nu_2|)\right)^{-1/2} + O\left( ((1+|\nu_0|) (1+|\nu_1|) (1+|\nu_2|))^{-3/4}\right)\\
& \gg \left((1+|\nu_0|) (1+|\nu_1|) (1+|\nu_2|)\right)^{-1/2}  
 \end{split}
\end{displaymath}
if one of $\nu_1, \nu_2$ is sufficiently large. \qed \\ 
 
 We briefly discuss  cusp forms $u : SL_2(\Bbb{Z}) \backslash \mathfrak{h}^2\rightarrow \Bbb{C}$ for the group $SL_2(\Bbb{Z})$ and spectral parameter $\nu \in i \Bbb{R}$ (Selberg's eigenvalue conjecture is known for $SL_2(\Bbb{Z})$).  A cusp form $u$ has  a Fourier expansion
 \begin{displaymath}
   u(z) = \sum_{m \not= 0} \frac{\rho_u(m)}{\sqrt{m}} W_{\nu}(|m|y) e(mx)
 \end{displaymath}
 where $W_{\nu}$ was defined in \eqref{gl2}. 
Similarly as in Lemma \ref{lem1} we see that an arithmetically normalized newform $u$ has norm
\begin{equation}\label{normgl2}
  \| u \|^2 = \int_{SL_2(\Bbb{Z}) \backslash \mathfrak{h}^2}  |u(z)|^2 \frac{dx\, dy}{y^2} = 2 L(\text{\rm Ad}^2 u, 1). 
\end{equation}
Indeed, the Eisenstein series $E(z, s) = \frac{1}{2}\displaystyle \sum_{\gamma \in  U_2\backslash  SL_2(\Bbb{Z})} \Im(\gamma z)^s$ has residue $3/\pi$ at $s=1$, hence by \eqref{gl2}
\begin{displaymath}
\begin{split}
  \| u \|^2 &= \frac{\pi}{3} \underset{s=1}{\text{res}}\sum_{m \not= 0} \frac{|\rho(m)|^2}{|m|} \frac{4 \pi}{ |\Gamma(1/2 + \nu)|^{2}} \int_{-\infty}^{\infty} |m|y K_{\nu}(2 \pi |m|y) K_{\bar{\nu}}(2 \pi |m|y) y^s \frac{dy}{y^2}\\
&  =  \frac{2\pi}{3 \zeta(2)} L({\rm Ad}^2 u, 1)  \frac{4 \pi }{|\Gamma(1/2 + \nu)|^2} \frac{\Gamma(\frac{1+\nu+\bar{\nu}}{2}) \Gamma(\frac{1-\nu+\bar{\nu}}{2})\Gamma(\frac{1+\nu-\bar{\nu}}{2})\Gamma(\frac{1-\nu-\bar{\nu}}{2})}{8}  = 2L(\text{Ad}^2u, 1);
\end{split}
\end{displaymath}
the evaluation of the integral follows from \cite[6.576.4]{GR} or Stade's formula for $GL(2)$. Again we see that an arithmetically normalized cusp form $u$ is essentially $L^2$-normalized, and $W_{\nu}$ has roughly norm one with respect to the inner product
 \begin{displaymath}
  (f, g) := \int_0^{\infty} f(y) \overline{g(y)} \det\left(\begin{smallmatrix} y  & \\   & 1\end{smallmatrix}\right)  \frac{dy}{y^2}.
 \end{displaymath}
 
\section{Eisenstein series}

There are three types of Eisenstein series on the space $L^2(\Gamma\backslash \mathfrak{h}^3)$ according to the decomposition
\begin{displaymath}
  L^2(SL_2(\Bbb{Z})\backslash \mathfrak{h}^2) =  L^2_{\text{Eis}} \oplus L^2_{\text{cusp}}  \oplus \Bbb{C} \cdot \textbf{1}. 
\end{displaymath}

The first  term gives rise to \emph{minimal parabolic Eisenstein series}: for $z \in \mathfrak{h}^3$ and $\Re \nu_1, \Re \nu_2$ sufficiently large we define the minimal parabolic Eisenstein series
\begin{displaymath}
  E(z, \nu_1, \nu_2) :=    \sum_{\gamma \in U_3 \backslash  \Gamma} I_{\nu_1, \nu_2}(\gamma z)
\end{displaymath}
where $I_{\nu_1, \nu_2}$ was defined in \eqref{defI}. It has meromorphic continuation in $\nu_1$ and $\nu_2$, and its non-zero Fourier coefficients are given by 
\begin{equation}\label{foureis1}
\begin{split}
  \int_0^1\int_0^1 \int_0^1 E(z, \nu_1, \nu_2) e(-m_1x_1-m_2x_2) dx_1 \, dx_2 \, dx_3& = \frac{A_{(\nu_1, \nu_2)}(m_1, m_2)}{|m_1m_2|}  \frac{W_{\nu_1, \nu_2}(m_1y_1, |m_2|y_2)}{\zeta(1+3\nu_0)\zeta(1+3\nu_1)\zeta(1+3\nu_2)}\\
  & = \frac{A_{(\nu_1, \nu_2)}(m_1, m_2)}{|m_1m_2|}  \frac{\tilde{W}_{\nu_1, \nu_2}(m_1y_1, |m_2|y_2)c_{\nu_1, \nu_2}^{-1}}{\zeta(1+3\nu_0)\zeta(1+3\nu_1)\zeta(1+3\nu_2)}
  \end{split}
\end{equation}
(cf.\ \eqref{defW} and \eqref{const} for the notation) where
\begin{displaymath}
  A_{(\nu_1, \nu_2)}(m_1, m_2)  = |m_1|^{\nu_1+2\nu_2} |m_2|^{2\nu_1+\nu_2} \sigma_{-3\nu_2, -3\nu_1}(|m_1|, |m_2|) 
\end{displaymath}
and $\sigma_{\nu_1, \nu_2}(m_1, m_2)$ is the multiplicative function defined by
\begin{displaymath}
  \sigma_{\nu_1, \nu_2}(p^{k_1}, p^{k_2}) = p^{-\nu_2k_1} \frac{\left|\left(\begin{smallmatrix} 1 & p^{\nu_2(k_1+k_2+2)} & p^{(\nu_1+\nu_2)(k_1+k_2 + 2}\\ 1 & p^{\nu_2(k_1+1)}  & p^{(\nu_1+\nu_2)(k_1+1)}\\ 1 & 1 & 1 \end{smallmatrix}\right)\right|}{\left|\left(\begin{smallmatrix} 1 & p^{2\nu_2} & p^{2(\nu_1+\nu_2)}\\ 1 & p^{\nu_2 }  & p^{\nu_1+\nu_2 }\\ 1 & 1 & 1 \end{smallmatrix}\right)\right|}.
\end{displaymath}
This is a combination of \cite[(6.5), (6.7), (6.8), (7.3), Theorem 7.2]{Bu}.  An alternative description is given as follows: $A_{(\nu_1, \nu_2)}(m_1, m_2)$ is defined by
\begin{displaymath}
  A_{(\nu_1, \nu_2)}(1, m) = \sum_{d_1d_2d_3 = m} d_1^{\alpha_1} d_2^{\alpha_2} d_3^{\alpha_3}
  \end{displaymath}
and the symmetry and Hecke relation
\begin{equation}\label{hecke}
\begin{split}
& A_{(\nu_1, \nu_2)}(m, 1) = \overline{A_{(\nu_1, \nu_2)}(1, m)} = A_{(\nu_2, \nu_1)}(1, m),\\
&  A_{(\nu_1, \nu_2)}(n, m) = \sum_{d \mid (n, m) } \mu(d) A_{(\nu_1, \nu_2)}\left(\frac{n}{d}, 1\right)A_{(\nu_1, \nu_2)}\left(1, \frac{m}{d}\right),
  \end{split}
\end{equation}
cf.\ \cite[Theorem 6.4.11 and Theorem 10.8.6]{Go} and note his different choice of the $I$-function. \\

Next we define \emph{maximal parabolic Eisenstein series}. Let $s \in \Bbb{C}$ have sufficiently large real part and let $u : SL_2(\Bbb{Z}) \backslash \mathfrak{h}^2 \rightarrow\Bbb{C}$ be a Hecke-Maa{\ss} cusp form with  $\|u \| = 1$,  spectral parameter  $ \nu \in  i \Bbb{R}$   and Hecke eigenvalues $\lambda_u(m)$.    Then we define
\begin{equation}\label{twistedeis}
  E(z,s; u) := \sum_{\gamma \in  P\backslash \Gamma} \det(\gamma z)^{s} u(\pi(\gamma z))
\end{equation}
where 
\begin{displaymath}
  \pi : \mathfrak{h^3} \rightarrow \mathfrak{h}^2, \quad  \left(\begin{matrix} y_1y_2 & x_2y_1 & x_3\\ & y_1 & x_1\\ & & 1\end{matrix} \right)  \mapsto \left(\begin{matrix} y_2 & x_2  \\ & 1\end{matrix} \right)
\end{displaymath}
is the restriction to the upper left corner. It has a meromorphic continuation in $s$, and as the minimal parabolic Eisenstein series it is an eigenform of all Hecke operators; in particular for $s = 1/2 + \mu$ it is an eigenform of $T(1, m)$ with eigenvalue
\begin{displaymath}
    B_{(\mu, u)}(1, m) = \sum_{d_1d_2 = |m|} \lambda_u(d_1) d_1^{-\mu}d_2^{2\mu}, 
\end{displaymath}
see \cite[Proposition 10.9.3]{Go}. We extend this definition to all pairs of integers by the Hecke relations \eqref{hecke}. Coupling this with \cite[Proposition 10.9.1]{Go}, we conclude that the non-zero Fourier coefficients  
\begin{displaymath}
  \int_0^1\int_0^1 \int_0^1 E(z, 1/2 + \mu; u) e(-m_1x_1-m_2x_2) dx_1 \, dx_2 \, dx_3
\end{displaymath}
are proportional to
\begin{equation}\label{foureis2}  
    \frac{B_{(\mu, u)}(m_1, m_2)}{|m_1m_2|}  \tilde{W}_{\mu - \frac{1}{3}\nu, \frac{2}{3} \nu}(m_1y_1, |m_2|y_2), 
\end{equation}
and the proportionality constant is
\begin{equation}\label{foureis2a}
\frac{c }{L(u, 1 + 3 \mu) L( \text{Ad}^2 u, 1)^{1/2} } 
\end{equation}  
  for some absolute non-zero  constant $c$. This can be seen by setting $m_1 = m_2 = 1$ and comparing with \cite[Theorem 7.1.2]{Sh} in the special case $G = GL(3)$, $M = GL(2) \times GL(1)$, $m=1$, $s = 3\mu$ and observing \eqref{normgl2}. \\

A degenerate case of \eqref{twistedeis} occurs if we choose $\phi$ to be the constant function (of course, this is not a cusp form). For $\Re s$ sufficiently large and $z \in \mathfrak{h^3}$ let 
\begin{equation}\label{eisdegen}
  E(z, s, \textbf{1}) := \sum_{\gamma \in P\backslash  \Gamma} \det(\gamma z)^s.  
\end{equation}
 This function has a meromorphic continuation to all  $s \in \Bbb{C}$, and it has a simple pole at $s=1$ with constant residue
\begin{equation}\label{res}
  \underset{s=1}{\text{res}} E(z, s, \textbf{1}) = \frac{1}{3} \left( \Gamma(3/2)\zeta(3)\pi^{-3/2}\right)^{-1}  
  = \frac{2 \pi }{3\zeta(3)},  
  \end{equation}
see \cite[Corollary 2.5]{Fr}\footnote{Note that the Eisenstein series in \cite[p.\ 164]{Fr} differs by a factor two from our definition. In \cite[Theorem 7.4.4]{Go} our definition is used, but the factor 1/2 seems to have got lost in the last display of p.\ 224 and the following argument.}. As the constant function on $SL_2(\Bbb{Z}) \backslash \mathfrak{h}^2$ is the residue of an Eisenstein series on $SL_2(\Bbb{Z}) \backslash \mathfrak{h}^2$, the Eisenstein series \eqref{eisdegen}  is a residue of a minimal parabolic Eisenstein series and  has only degenerate terms in its Fourier expansion. \\

\section{Kloosterman sums}\label{Kloosterman}

As usual we write
\begin{displaymath}
  S(m, n, c) := \underset{d \, (\text{mod }c)}{\left.\sum\right.^{\ast}}e\left(\frac{md + n\bar{d}}{c}\right)
\end{displaymath}
for the standard Kloosterman sum. We introduce now $GL(3)$ Kloosterman sums; the following account is taken from \cite{BFG}. \\

For $n_1, n_2, m_1, m_2 \in \Bbb{Z}$, $D_1, D_2 \in \Bbb{N}$ we define
\begin{equation}\label{def}
\begin{split}
 & S(m_1, m_2, n_1, n_2, D_1, D_2) := \\
  &\underset{\substack{B_1, C_1 \, (\text{mod }D_1)\\ B_2, C_2 \, (\text{mod }D_2)\\ (D_1, B_1, C_1) = (D_2, B_2, C_2) = 1\\ D_1C_2 + B_1B_2 + C_1D_2 \equiv 0 \, (D_1D_2)}}{\sum\sum\sum\sum} e\left(\frac{m_1B_1 + n_1(Y_1D_2 - Z_1 B_2)}{D_1} \right) e\left(\frac{m_2B_2 + n_2(Y_2 D_1 - Z_2 B_1)}{D_2}\right)
  \end{split}
\end{equation}
where $Y_1, Y_2, Z_1, Z_2$ are chosen such that
\begin{displaymath}
  Y_1B_1 + Z_1C_1 \equiv 1 \, (\text{mod }D_1), \quad  Y_2B_2 + Z_2C_2 \equiv 1 \, (\text{mod } D_2).
\end{displaymath}
It can be shown that this expression is well-defined \cite[Lemma 4.1, 4.2]{BFG}. Clearly it depends only on $m_1, n_1 \, (\text{mod }D_1)$ and $m_2, n_2 \, (\text{mod }D_2)$, and satisfies \cite[Property 4.4, 4.5]{BFG}
\begin{equation}\label{property}
  S(m_1, m_2, n_1, n_2, D_1, D_2) = S(m_2, m_1, n_2, n_1, D_2, D_1) = S(n_1, n_2, m_1, m_2, D_1, D_2). 
\end{equation}
Moreover, if $p_1p_2 \equiv q_1q_2 \equiv 1$ (mod $D_1 D_2$), then \cite[Property 4.3]{BFG}
\begin{displaymath} 
  S(p_1m_1, p_2m_2, q_1m_1, q_2m_2, D_1, D_2) = S(m_1, m_2, m_1, m_2, D_1, D_2). 
\end{displaymath}
Finally we have the factorization rule \cite[Property 4.7]{BFG}
\begin{equation}\label{mult}
\begin{split}
 & S(m_1, m_2, n_1, n_2, D_1D_1', D_2D_2') \\
  & = S(\overline{D'_1}^2 D_2'm_1, \overline{D_2'}^2D_1'm_2, n_1, n_2, D_1, D_2) S(\overline{D_1}^2 D_2m_1, \overline{D_2}^2D_1m_2, n_1, n_2, D_1', D_2') 
  \end{split}
\end{equation}
whenever $(D_1D_2, D_1'D_2') = 1$ and inverses are taken with respect to the product of the respective moduli, that is,
\begin{displaymath}
  \overline{D_1} D_1 \equiv \overline{D_2}D_2 \equiv 1 \, (\text{mod }D_1'D_2'), \quad \overline{D_1'} D_1' \equiv \overline{D_2'}D_2' \equiv 1 \, (\text{mod }D_1D_2).
\end{displaymath}
This implies in particular 
\begin{equation}\label{decomp}
  S(m_1, m_2, n_1, n_2, D_1, D_2) = S(D_2m_1, n_1, D_1)S(D_1 m_2, n_2, D_2), \quad (D_1, D_2) = 1. 
\end{equation}
For a prime $p$ and $l \geq 1$ we have \cite[Property 4.10]{BFG}
\begin{equation}\label{p}
  S(m_1, m_2, n_1, n_2, p, p^{l}) = S(n_1, 0, p)S(m_2, n_2p, p^l) + S(m_1, 0, p)S(n_2, m_2p, p^l) + \delta_{l=1}(p-1). 
\end{equation}
Essentially best possible (``Weil-type") upper bounds for $S(m_1, m_2, n_1, n_2, D_1, D_2)$ have been given by Stevens \cite[Theorem 5.1]{Ste}. The dependence on $m_1, m_2, n_1, n_2$ has been worked out in \cite[p.\ 39]{Buth}.

\begin{lemma}\label{prop4} For any integers $n_1, n_2, m_1, m_2 \in \Bbb{Z} \setminus \{0\}$, $D_1, D_2 \in \Bbb{N}$ and any $\varepsilon > 0$  we have
\begin{equation*} 
S(m_1, m_2, n_1, n_2, D_1, D_2) \ll  (D_1D_2)^{1/2+\varepsilon} \bigl((D_1, D_2) (m_1n_2, [D_1, D_2])(m_2n_1, [D_1, D_2])\bigr)^{1/2}
\end{equation*}
where $[., .]$ denotes the least common multiple. In particular,
 \begin{equation}\label{kloosum}
   \sum_{D_1 \leq X_1} \sum_{D_2 \leq X_2} |S(m, \pm 1, n, 1, D_1, D_2)|  \ll (X_1X_2)^{3/2+\varepsilon} (n, m)^{\varepsilon}
 \end{equation}
 if $mn \not= 0$. All implied constants depend only on $\varepsilon$.
\end{lemma}

\textbf{Proof.} It remains to show \eqref{kloosum} which is straightforward:
\begin{displaymath}
\begin{split}
   \sum_{D_1 \leq X_1} \sum_{D_2 \leq X_2} & |S(m, \pm 1, n, 1, D_1, D_2)|  \ll  (X_1X_2)^{1/2+\varepsilon} \sum_{d_1 \mid m} \sum_{d_2 \mid n} (d_1d_2)^{1/2} \underset{[d_1, d_2] \mid DD_1D_2}{\sum_{D} \sum_{D_1 \leq X_1/D} \sum_{D_2 \leq X_2/D_2}} D^{1/2} \\
   & \ll (X_1X_2)^{3/2+\varepsilon} \sum_{d_1 \mid m} \sum_{d_2 \mid n} \frac{(d_1d_2)^{1/2} }{[d_1, d_2]} \ll (X_1X_2)^{3/2+\varepsilon} (n, m)^{\varepsilon}. 
 \end{split}  
\end{displaymath}
  \qed \\
 
Next we define a different class of Kloosterman sums: If $D_1 \mid D_2$, we put
\begin{displaymath}
  \tilde{S}(m_1, n_1, n_2, D_1, D_2) := \underset{\substack{ C_1 \, (D_1), \,\, C_2 \, (D_2)\\(C_1, D_1) = (C_2, D_2/D_1) = 1}}{\sum\sum} e\left(\frac{m_1C_1 + n_1 \overline{C_1}C_2}{D_1}\right) e\left(\frac{n_2\overline{C_2}}{D_2/D_1}\right).
\end{displaymath}
Again this sum depends only on $m_1, n_1 \, (\text{mod }D_1)$ and $n_2 \, (\text{mod }D_2/D_1)$, and for $p_1q_1 \equiv 1 \, (\text{mod }D_1)$, $p_2q_2 \equiv 1 \, (\text{mod }D_2)$ we have \cite[Property 4.13]{BFG}
\begin{displaymath}
  \tilde{S}(m_1p_1, n_1q_1p_2, n_2q_2, D_1, D_2) = \tilde{S}(m_1, n_1, n_2, D_1, D_2).
\end{displaymath}
We have the factorization rule \cite[Property 4.15]{BFG}
\begin{displaymath}
  \tilde{S}(m_1, n_1, n_2, D_1 D_1', D_2 D_2') = \tilde{S}(m_1 \overline{D_1'}, n_1D_2', n_2 \overline{D_2'}^2, D_1, D_2) \tilde{S}(m_1 \overline{D_1}, n_1 D_2, n_2\overline{D_2}^2, D_1', D_2')
\end{displaymath}
whenever $(D_2, D_2') = 1$ and all terms are defined. Finally  we have for a prime number $p$ and $1 \leq l < k$ \cite[Property 4.16, 4.17]{BFG}
\begin{displaymath}
 \tilde{S}(m_1, n_1, n_2, p^l, p^l) = \begin{cases} p^{2l} - p^{2l-1}, & p^l \mid m_1, \, p^l \mid n_1\\
 -p^{2l-1}, & p^{l-1} \parallel m_1, \, p^l \mid n_1\\
 0, &\text{otherwise}
 \end{cases}
\end{displaymath}
and 
\begin{displaymath}
   \tilde{S}(m_1, n_1, n_2, p^l, p^k) = 0
\end{displaymath}
unless
\begin{itemize}
\item $k < 2l$ and $p^{2l-k} \mid n_1$, or\\
\item $k = 2l$, or\\
\item $k > 2l$ and $p^{k-2l} \mid n_2$. 
\end{itemize}
In particular
\begin{displaymath}
  \tilde{S}(m_1, n_1, n_2, D_1, D_2) = 0 \quad \text{unless} \quad D_1^2 \mid n_1D_2. 
\end{displaymath}

A sharp bound was proved by Larsen \cite[Appendix]{BFG}:
\begin{equation}\label{sharp}
   \tilde{S}(m_1, n_1, n_2, D_1, D_2) \ll \min\left((n_2, D_2/D_1)D_1^2, (m_1, n_1, D_1)D_2\right)(D_1D_2)^{\varepsilon}. 
\end{equation}

\section{Poincar\'e series}

Let $F : (0, \infty)^2 \rightarrow \Bbb{C}$ be a smooth compactly supported function (or sufficiently rapidly decaying at 0 and $\infty$ in both variables). Let
\begin{displaymath}
  F^{\ast}(y_1, y_2) := F(y_2, y_1). 
\end{displaymath}
 For two positive integers $m_1, m_2$ and $z \in \mathfrak{h}^3$ let $\mathcal{F}_{m_1, m_2}(z) := e(m_1x_2 + m_2x_2) F(m_1y_1, m_2y_2)$. Then we consider   the following Poincar\'e series:
\begin{displaymath}
  P_{m_1, m_2}(z) := \sum_{\gamma \in  U_3 \backslash  \Gamma} \mathcal{F}_{m_1, m_2}(\gamma z). 
\end{displaymath}
Unfolding shows
\begin{equation}\label{unfold1}
\begin{split}
  \langle \phi, P_{m_1, m_2} \rangle & = \int_{U_3 \backslash \mathfrak{h}^3} \phi(z) \overline{\mathcal{F}_{m_1, m_2}(z)} \, dx_1\, dx_2\, dx_3 \frac{dy_1\, dy_2}{(y_1y_2)^3}\\
  &=  \int_0^{\infty} \int_0^{\infty}   \int_0^1\int_0^1 \int_0^1 \phi(z) e(-m_1x_1-m_2x_2) dx_1\, dx_2\, dx_3   \, \overline{ F(m_1y_1, m_2y_2) } \frac{dy_1\, dy_2}{(y_1y_2)^3} \\
  \end{split}
\end{equation}
for an arbitrary automorphic form $\phi$. In particular, if $\phi$ is given as in \eqref{defphi}, we find 
\begin{equation}\label{unfold2}
\begin{split}
  \langle \phi, P_{m_1, m_2} \rangle  &=  \int_0^{\infty} \int_0^{\infty}    \frac{A_{\phi}(m_1, m_2)}{m_1m_2} \tilde{W}_{\nu_1, \nu_2}(m_1y_1, m_2y_2) \overline{F(m_1y_1, m_2y_2)}  \frac{dy_1\, dy_2}{(y_1y_2)^3} \\
  & = m_1m_2 A_{\phi}(m_1, m_2)  \int_0^{\infty} \int_0^{\infty}   \tilde{W}_{\nu_1, \nu_2}(y_1, y_2) \overline{F(y_1, y_2) } \frac{dy_1\, dy_2}{(y_1y_2)^3}. 
  \end{split}
\end{equation}
We want to apply \eqref{unfold1} also with $\phi = P_{n_1, n_2}$ where $n_1, n_2$ is another pair of positive integers. The Fourier expansion of $P_{n_1, n_2}$ has been computed explicitly in \cite[Theorem 5.1]{BFG}: For $m_1, m_2 > 0$ we have
\begin{equation}\label{unfold3}
 \int_0^1\int_0^1 \int_0^1 P_{n_1, n_2}(z) e(-m_1x_1 - m_2x_2) dx_1 \, dx_2 \, dx_3  = S_1 + S_{2a} + S_{2b} + S_3,
\end{equation} 
where
\begin{displaymath}
\begin{split}
  S_1 &= \delta_{m_1, n_1} \delta_{m_2, n_2} F(n_1y_1, n_2 y_2),\\
   S_{2a}& = \sum_{\epsilon  = \pm 1} \sum_{\substack{  D_1 \mid D_2\\  m_2D_1^2=n_1D_2}}\tilde{S}(\epsilon m_1, n_1, n_2, D_1, D_2) \tilde{J}_F\left(y_1, y_2, \epsilon m_1, n_1, n_2, D_1, D_2 \right), \\
   S_{2b} &=   \sum_{\epsilon  = \pm 1} \sum_{\substack{  D_2 \mid D_1\\ m_1D_2^2 = n_2D_1}}\tilde{S}(\epsilon m_2, n_2, n_1, D_2, D_1) \tilde{J}_{F^{\ast}}\left(y_2, y_1,\epsilon m_2 ,n_2, n_1, D_2, D_1 \right), \\
   S_3 &= \sum_{\epsilon_1, \epsilon_2 = \pm 1} \sum_{D_1, D_2 } S(\epsilon_1m_1, \epsilon_2m_2, n_1, n_2, D_1, D_2) J(y_1, y_2, \epsilon_1m_1, \epsilon_2 m_2, n_1, n_2, D_1, D_2). 
  \end{split}
\end{displaymath}
The Kloosterman sums have been defined in Section \ref{Kloosterman} and the weight functions are given as follows:
\begin{equation}\label{prelimj1}
\begin{split}
  \tilde{J}_F&(y_1, y_2,  m_1, n_1, n_2, D_1, D_2) = y_1^2y_2 \int_{\Bbb{R}^2} e(-m_1x_1y_1) e\left(\frac{n_1D_2y_2}{D_1^2} \cdot \frac{x_1x_2}{x_1^2 + 1}\right)  \\
& \times e\left(\frac{n_2D_1}{y_1y_2 D_2^2} \cdot \frac{x_2}{x_1^2 + x_2^2 + 1}\right) F\left(\frac{n_1D_2y_2}{D_1^2} \cdot \frac{\sqrt{x_1^2+x_2^2 + 1}}{x_1^2 + 1}, \frac{n_2D_1}{y_1y_2 D_2^2}\cdot \frac{\sqrt{x_1^2 + 1}}{x_1^2+x_2^2 + 1}\right) dx_1\, dx_2, 
\end{split}  
  \end{equation}
\begin{equation}\label{prelimj2}
\begin{split}
 & J(y_1, y_2, m_1, m_2, n_1, n_2, D_1, D_2) = (y_1y_2)^2 \\
& \times \int_{\Bbb{R}^3} e\left(-m_1x_1y_1 - m_2x_2y_2\right) e\left(-\frac{n_1D_2}{D_1^2y_2} \cdot \frac{x_1x_3+x_2}{x_3^2 + x_2^2 + 1}\right) e\left(-\frac{n_2D_1}{D_2^2y_1} \cdot \frac{x_2(x_1x_2-x_3) + x_1} {(x_1x_2-x_3)^2 + x_1^2 + 1}\right) \\
& \quad\quad\quad \times F\left(\frac{n_1D_2}{D_1^2y_2} \cdot \frac{\sqrt{(x_1x_2-x_3)^2 + x_1^2 + 1}}{x_3^2 + x_2^2 + 1}, \frac{n_2D_1}{D_2^2y_1} \cdot \frac{ \sqrt{x_3^2 + x_2^2 + 1}} {(x_1x_2-x_3)^2 + x_1^2 + 1}\right) dx_1\, dx_2\, dx_3.
\end{split}  
\end{equation}

\section{Spectral decomposition}

We have the following spectral decomposition theorem \cite[Proposition 10.13.1]{Go}: If $\phi \in L^2(\Gamma \backslash \mathfrak{h}^3)$ is orthogonal to all residues of Eisenstein series, then
\begin{displaymath}
\begin{split}
 \phi  = & \sum_{j} \langle \phi, \phi_j\rangle \phi_j + \int_{(0)} \int_{(0)} \langle \phi, E(., \nu_1, \nu_2)\rangle  E(., \nu_1, \nu_2) \frac{d\nu_1\, d\nu_2}{(4\pi i)^2} \\
& +   \sum_{j} \int_{(0)} \langle \phi, E(., 1/2 + \mu;  u_j)\rangle  E(., 1/2 + \mu; u_j ) d\mu  + \frac{1}{2\pi i}   \int_{(0)} \langle \phi, E(., 1/2 + \mu;  \textbf{1})\rangle  E(.,1/2 + \mu; \textbf{1} ) \frac{d\mu}{2\pi i}
 \end{split}
\end{displaymath}
where the first $j$-sum runs over an orthonormal basis of cusp forms $\phi_j$ for $SL_3(\Bbb{Z})$ and the second $j$-sum runs over an orthonormal basis of cusp forms $u_j$ for $SL_2(\Bbb{Z})$. 

Therefore we have for 4 positive integers $n_1, n_2, m_1, m_2$ an equality of the type
\begin{displaymath}
 \frac{ \langle P_{n_1, n_2},  P_{m_1, m_2} \rangle }{n_1n_2m_1m_2}=   \sum_{j}  \frac{\langle P_{n_1, n_2}, \phi_j\rangle \langle \phi_j, P_{m_1, m_2} \rangle}{n_1n_2m_1m_2 
 } + \ldots \text{(continuous spectrum)}.
\end{displaymath}
We refer to the right hand side as the spectral side and to the left hand side as the arithmetic side. 

We proceed to describe the spectral side and the arithmetic side more precisely. We define an inner product on $L^2((0, \infty)^2, dy_1 dy_2/(y_1y_2)^3)$ by
\begin{displaymath}
  \langle f, g \rangle := \int_0^{\infty} \int_0^{\infty}  f(y_1, y_2) \overline{g(y_1, y_2)} \frac{dy_1 dy_2}{(y_1y_2)^3}. 
\end{displaymath}
Let $\{\phi_j\}$ denote an arithmetically normalized ortho\emph{gonal} basis of the space of cusp forms on $L^2(SL_3(\Bbb{Z}) \backslash \mathfrak{h}^3)$ that we assume to be eigenfunctions of the Hecke algebra with eigenvalues $A_j(m_1, m_2)$. Let $\{u_j\}$ be an  arithmetically normalized orthogonal basis of the space of cusp forms on $L^2(SL_2(\Bbb{Z})\backslash \mathfrak{h}^2)$ that we  assume to be eigenfunctions of the Hecke algebra with eigenvalues $\lambda_j(m)$ and spectral parameter $\nu_j \in i \Bbb{R}$.  

\begin{prop}\label{kuznetsov} Keep the notation developed so far.  Let $F : (0, \infty)^2 \rightarrow \Bbb{C}$ be a smooth compactly supported function, and let $m_1, m_2, n_1, n_2 \in \Bbb{N}$. Then for some absolute constant $c > 0$ the following equality holds:
\begin{equation}\label{spectral}
\begin{split}
   & \sum_{j}  \frac{\overline{A_{j}(n_1, n_2)} A_{j}(m_1,m_2)}{\| \phi_j \|^2}    | \langle \tilde{W}_{\nu_1, \nu_2}, F \rangle |^2 \\
    &+ \frac{1}{(4\pi i)^2} \int_{(0)} \int_{(0)}  \frac{\overline{A_{(\nu_1, \nu_2)}(n_1, n_2)} A_{(\nu_1, \nu_2)}(m_1,m_2)}{|\zeta(1+3\nu_0)\zeta(1+3\nu_1)\zeta(1+3\nu_2)|^2}    | \langle \tilde{W}_{\nu_1, \nu_2}, F \rangle |^2 d\nu_1 d\nu_2\\
 & +    \frac{c}{2\pi i} \sum_{j} \int_{(0)}     \frac{\overline{B_{(\mu, u_j)}(n_1, n_2)} B_{(\mu, u_j)}(m_1,m_2)}{|L(u_j, 1+3\mu)|^2 L(\text{\rm Ad}^2 u_j, 1)\ }    | \langle \tilde{W}_{\mu - \frac{1}{3}\nu_j, \frac{2}{3}\nu_j}, F \rangle |^2
  d\mu \\
 & =  \Sigma_1 + \Sigma_{2a} + \Sigma_{2b} + \Sigma_3, 
    \end{split}
\end{equation}
 where
\begin{equation}\label{arithmetic}
\begin{split}
  \Sigma_1& = \delta_{m_1, n_1} \delta_{m_2, n_2}   \| F \|^2,\\
   \Sigma_{2a}& =  \sum_{\epsilon  = \pm 1} \sum_{\substack{  D_1 \mid D_2\\  m_2D_1^2= n_1D_2}}\frac{\tilde{S}(\epsilon m_1, n_1, n_2, D_1, D_2)}{D_1D_2}  \tilde{\mathcal{J}}_{\epsilon; F}\left(\sqrt{\frac{n_1n_2m_1}{D_1D_2}} \right),\\
 \Sigma_{2b} &=    \sum_{\epsilon  = \pm 1} \sum_{\substack{  D_2 \mid D_1\\ m_1D_2^2 = n_2D_1}} \frac{\tilde{S}(\epsilon m_2, n_2, n_1, D_2, D_1) }{D_1D_2}\tilde{\mathcal{J}}_{\epsilon; F^{\ast}}\left( \sqrt{\frac{n_1n_2m_2}{D_1D_2}}  \right), \\
   \Sigma_3 &=  \sum_{\epsilon_1, \epsilon_2 = \pm 1} \sum_{D_1, D_2  } \frac{S(\epsilon_1m_1, \epsilon_2m_2, n_1, n_2, D_1, D_2)}{D_1D_2}\mathcal{J}_{\epsilon_1, \epsilon_2}\left( \frac{\sqrt{m_1n_2 D_1}}{D_2 }, \frac{\sqrt{m_2n_1D_2}}{D_1}\right).
  \end{split}
\end{equation}
The weight functions $\tilde{\mathcal{J}}$ and $\mathcal{J}$ are given by
\begin{equation}\label{test1}
\begin{split}
 & \tilde{\mathcal{J}}_{\epsilon;F}(A) =  A^{-2} \int_0^{\infty}\int_0^{\infty} \int_{-\infty}^{\infty}\int_{-\infty}^{\infty}   e(-\epsilon Ax_1y_1) e\left(y_2 \cdot \frac{x_1x_2}{x_1^2 + 1}\right) e\left(\frac{A}{y_1y_2} \cdot \frac{x_2}{x_1^2 + x_2^2 + 1}\right) \\
& \quad\quad\quad F\left(y_2 \cdot \frac{\sqrt{x_1^2+x_2^2 + 1}}{x_1^2 + 1}, \frac{A}{y_1y_2}\cdot \frac{\sqrt{x_1^2 + 1}}{x_1^2+x_2^2 + 1}\right) \overline{F(Ay_1, y_2)} dx_1\, dx_2 \frac{dy_1\, dy_2}{y_1y_2^2},
\end{split}  
\end{equation}
\begin{equation}\label{test2}
\begin{split}
&  \mathcal{J}_{\epsilon_1, \epsilon_2} (A_1, A_2) =  (A_1A_2)^{-2} \int_0^{\infty}\int_0^{\infty} \int_{-\infty}^{\infty}\int_{-\infty}^{\infty}\int_{-\infty}^{\infty} e\left(-\epsilon_1 A_1x_1y_1 - \epsilon_2 A_2x_2y_2\right)\\
 & \times  e\left(-\frac{A_2}{y_2} \cdot \frac{x_1x_3+x_2}{x_3^2 + x_2^2 + 1}\right) e\left(-\frac{A_1}{y_1} \cdot \frac{x_2(x_1x_2-x_3) + x_1} {(x_1x_2-x_3)^2 + x_1^2 + 1}\right)\overline{F(A_1y_1, A_2 y_2)} \\
& \times F\left(\frac{A_2}{y_2} \cdot \frac{\sqrt{(x_1x_2-x_3)^2 + x_1^2 + 1}}{x_3^2 + x_2^2 + 1}, \frac{A_1}{y_1} \cdot \frac{ \sqrt{x_3^2 + x_2^2 + 1}} {(x_1x_2-x_3)^2 + x_1^2 + 1}\right)  dx_1\, dx_2\, dx_3 \frac{dy_1\, dy_2}{y_1y_2}.
\end{split}  
\end{equation}
\end{prop}

\textbf{Proof.} The spectral side \eqref{spectral} follows from\footnote{Even though $\tilde{W}_{\nu_1, \nu_2}$ just fails to be in $L^2((0, \infty)^2, dy_1 dy_2/(y_1y_2)^3)$, the inner products exist by the decay properties of $F$.}  \eqref{unfold2} in combination with \eqref{foureis1} and \eqref{foureis2}, \eqref{foureis2a}. Note that $E(z, 1/2 + \mu, \textbf{1})$ does not contribute because it has only degenerate terms in its Fourier expansion. 

Upon combining \eqref{unfold1} and \eqref{unfold3}, we obtain the arithmetic side \eqref{arithmetic} after applying a linear change of variables 
\begin{displaymath}
   y_1 \mapsto \sqrt{\frac{n_1n_2}{D_1D_2m_1}} y_1, \quad y_2 \mapsto \frac{y_2}{m_2}
\end{displaymath}
in     \eqref{prelimj1} and observing $m_2D_1^2 = n_1D_2$ (and with interchanged indices for $\Sigma_{2b}$),  and
\begin{displaymath}
  y_1 \mapsto \sqrt{\frac{n_2D_1}{m_1D_2^2}} y_1, \quad  y_2 \mapsto \sqrt{\frac{n_1D_2}{m_2D_1^2}} y_2
\end{displaymath}
in \eqref{prelimj2}. \qed\\

Formally \eqref{spectral} and \eqref{arithmetic}  resemble the $GL(2)$ Kuznetsov formula, but in its present form it is relatively useless as long as we do not understand the transforms $|\langle \tilde{W}_{\nu_1, \nu_2}, F\rangle|^2$ and $\mathcal{\tilde{J}}, \mathcal{J}$ for a given test function $F$. The present formulation has the important advantage that the weight functions on the arithmetic side \eqref{arithmetic} do not depend on $n_1, n_2, m_1, m_2, D_1, D_2$ individually, but only in a coupled fashion.  This is, of course, a well-known phenomenon in the $GL(2)$ world. 

We choose now  
\begin{equation}\label{defE}
 F(y_1, y_2) = F_{f, X_1, X_2, R_1, R_2, \tau_1, \tau_2}(y_1, y_2) := (R_1R_2(R_1+R_2))^{1/2} f(X_1y_1)f(X_2y_2)  y_1^{i(\tau_1 + 2\tau_2)}y_2^{i(\tau_2 + 2\tau_1)}
\end{equation}
for $X_1, X_2, R_1, R_2 \geq  1$, $\tau_1, \tau_2 \geq 0$, $\tau_1 + \tau_2 \geq 1$ and $f$ a fixed  smooth, nonzero, non-negative function with support in $[1, 2]$.  Analytic properties of $\langle \tilde{W}_{\nu_1, \nu_2}, F  \rangle$ have been obtained in Proposition \ref{whitint}.  We summarize   some bounds for  the weight functions occurring on the   arithmetic side in the following proposition.
\begin{prop}\label{prop6}
 With the notation developed so far, we have 
  \begin{equation}\label{diag}
  \|F \|^2 \asymp  (X_1X_2)^2 R_1 R_2 (R_1 + R_2).
\end{equation}
Let $C_1, C_2\geq 0$, $\varepsilon > 0$. Then 
\begin{equation}\label{J1}
  \tilde{\mathcal{J}}_{\epsilon; F}(A) \begin{cases} \ll X_1^2X_2 R_1R_2(R_1+R_2)\left(\frac{1+A^{2/3} }{\tau_1+\tau_2}\right)^{C_1}, \\
  = 0,  \quad \text{ if } A\leq (100X_1)^{-3/2} + (100 X_1X_2)^{-3/4}, \end{cases}
\end{equation}
and
\begin{equation}\label{J2}
\begin{split}
   \mathcal{J}_{\epsilon_1, \epsilon_2}&(A_1,  A_2)   \ll (X_1X_2)^{2}R_1R_2(R_1+R_2) \\
   &  \times \left(\frac{1+A_1^{4/3}A_2^{2/3}(X_1+X_2)}{\tau_1+\tau_2}\right)^{C_1} \left(\frac{1+A_2^{4/3}A_1^{2/3}(X_1+X_2)}{\tau_1+\tau_2}\right)^{C_2}  ((X_1+A_1)(X_2+A_2))^{\varepsilon},\\
     \mathcal{J}_{\epsilon_1, \epsilon_2}&(A_1,  A_2)   = 0, \quad \text{ if } \min(A_1A_2^2, A_2A_1^2) \leq (100 X_1X_2)^{-3/2}. 
   \end{split}
\end{equation}
In the special case when $A_1, A_2 \leq 1$, $X_1 = 1$, $X_2 = X \geq 1$, $R_1 + R_2 \asymp \tau_1 + \tau_2$ this can be improved to
\begin{equation}\label{improvedJ2}
   \mathcal{J}_{\epsilon_1, \epsilon_2}(A_1,  A_2)   \ll X^2 R_1R_2   \left(\frac{1+A_1^{4/3}A_2^{2/3}X}{\tau_1+\tau_2}\right)^{C_1} \left(\frac{1+A_2^{4/3}A_1^{2/3}X}{\tau_1+\tau_2}\right)^{C_2}  ((R_1+R_2)X)^{\varepsilon}.
\end{equation}
Let $g$ be a fixed smooth function with compact support in $(0, \infty)$.   Then for  $R_1, R_2 \gg 1$ sufficiently large   we have
 \begin{equation}\label{intJ2}
\begin{split}
 & \int_{0}^{\infty}\int_{0}^{\infty} g\left(\frac{\tau_1}{R_1}\right)g\left(\frac{\tau_2}{R_2}\right)  \mathcal{J}_{\epsilon_1, \epsilon_2}(A_1,  A_2) d\tau_1 d\tau_2\\ 
 & \ll (X_1X_2)^2R_1R_2(R_1+R_2)  \left(\frac{1+A_1^{4/3}A_2^{2/3}(X_1+X_2)}{R_1+R_2}\right)^{C_1} \left(\frac{1+A_2^{4/3}A_1^{2/3}(X_1+X_2)}{R_1+R_2}\right)^{C_2} \\
 &\times   (R_1R_2(X_1+A_1)(X_2+A_2))^{\varepsilon}.
  \end{split}
\end{equation}
On the left hand side we have suppressed the dependence of $\mathcal{J}_{\epsilon_1, \epsilon_2}$ on $\tau_1, \tau_2$. 
 \end{prop}
  
 \textbf{Remark 4.}    The bounds \eqref{J1}, \eqref{J2},  \eqref{intJ2} are  not best possible, but \eqref{improvedJ2} is likely to be best possible.  The important feature is that \eqref{J1} and \eqref{J2}  effectively bound $A_1, A_2$ from below, and therefore $D_1, D_2$ in \eqref{arithmetic} from above. For example, for the contribution of the long Weyl element, we can essentially assume
 \begin{displaymath}
   D_1 \leq \frac{(m_1n_2)^{1/3}(m_2n_1)^{2/3}}{\tau_1 + \tau_2}, \quad D_2 \leq \frac{(m_1n_2)^{2/3}(m_2n_1)^{1/3}}{\tau_1 + \tau_2}
 \end{displaymath}
if $X_1 = X_2 = 1$.  It is instructive to compare this with the $GL(2)$ situation: one can construct a sufficiently nice test function $h$ on the spectral side with  essential support on $[T, T+1]$   such that the integral transforms $h^{\pm}$ in \eqref{proto} are negligible unless
    $c \leq \frac{(nm)^{1/2}}{T}$. 
   
   The bound \eqref{intJ2} shows that integration over $\tau_1, \tau_2$ can be performed at almost no cost, in other words, we save a factor $(R_1R_2)^{1-\varepsilon}$ compared to trivial integration. \\
  
\textbf{Remark 5.} Choosing $f(y_1, y_2) = e^{-2\pi (y_1+y_2)} (y_1y_2)^{100} $ (say), the two $y$-integrals in \eqref{test2} can be computed explicitly using \cite[3.471.9]{GR}, giving two Bessel-$K$-functions with  general complex arguments. It is not clear how to take advantage of this fact. \\

 \textbf{Proof.} Equation  \eqref{diag} is clear. We proceed to prove \eqref{J1}. Let us write
 \begin{displaymath}
   \xi_1 := x_1^2 + 1, \quad \xi_2 = x_1^2 + x_2^2 + 1
 \end{displaymath}
in \eqref{test1}.  The support of $f$ restricts the variables to
 \begin{equation}\label{trunc}
 \begin{split}
  & (X_1A)^{-1}   \leq  y_1 \leq 2(X_1A)^{-1} , \quad  X_2^{-1} \leq y_2 \leq 2X_2^{-1}, \\
  & X_2/(2X_1) \leq \xi_2^{1/2} \xi_1^{-1} \leq 2X_2/X_1, \quad (A^2X_1X_2^2)^{-1}  \leq \xi_1^{1/2} \xi_2^{-1} \leq  8(A^2X_1X_2^2)^{-1}. 
 \end{split} 
 \end{equation}
 The second set of conditions in \eqref{trunc} implies $  \xi_1 \asymp  A^{4/3} X_1^2$ and  $\xi_2 \asymp A^{8/3}(X_1X_2)^2$. Hence the second part of \eqref{J1} is clear and a trivial estimation shows
 \begin{displaymath}
   \tilde{\mathcal{J}}_{\epsilon; F}(A) \ll R_1R_2(R_1+R_2) X_1^2X_2. 
 \end{displaymath}
In certain ranges this can be improved by  partial integration.   We have 
 \begin{displaymath}
 \begin{split}
 & \tilde{\mathcal{J}}_{\epsilon; F}(A) = \frac{R_1R_2(R_1+R_2) }{A^{ 2- i(\tau_1 - \tau_2)}}  \int_0^{\infty}\int_0^{\infty} \int_{-\infty}^{\infty}\int_{-\infty}^{\infty}   (x_1^2+1)^{  - \frac{3}{2}i\tau_2}(x_1^2 + x_2^2+1)^{ - \frac{3}{2}i\tau_1}  e\left( -Ay_1\epsilon x_1\right)  \\
 &\quad\quad\quad \times e\left(\frac{y_2x_1x_2 }{x_1^2 + 1}  \right)e\left(\frac{Ax_2}{y_1y_2(x_1^2 + x_2^2 + 1)}\right) \overline{f}(X_1Ay_1)   \overline{f}(X_2y_2)f\left(X_1y_2\frac{ \sqrt{x_1^2 + x_2^2 + 1}}{x_1^2 + 1}\right) 
 \\
& \quad\quad\quad \times  f\left(\frac{AX_2}{y_1y_2} \cdot \frac{x_2 + i \sqrt{x_1^2+1}}{x_1^2 + x_2^2 + 1}\right)  y_1^{  - 3i (\tau_1 + \tau_2)} y_2^{   - 3i\tau_1}dx_1\, dx_2 \frac{dy_1\, dy_2}{y_1y_2^2}.
\end{split}  
 \end{displaymath}
   We can assume that $C_1$ is an integer.  Then $C_1$ successive integrations by parts with respect to $y_1$ yield an additional factor 
 \begin{equation}\label{paren}
 \ll_{C_1} \left(\left(\frac{y_1}{\tau_1+\tau_2}\right)\left( A|x_1| + \frac{A|x_2|}{\xi_2y_1^2y_2} + \frac{1}{y_1}\right)\right)^{C_1} \ll_{C_1} \left(\frac{1+A^{2/3}}{\tau_1+\tau_2}\right)^{C_1}
\end{equation} 
in the support of $f$.  \\
 
The bound \eqref{J2} can be shown similarly, but the estimations are a little more involved. Here we write 
\begin{equation}\label{xi}
   \xi_1 := (x_1x_2  - x_3)^2 + x_1^2 + 1, \quad \xi_2 = x_3^2 + x_2^2 + 1
\end{equation}
and truncate
\begin{equation}\label{sizey}
  \begin{split}
  & (X_1A_1)^{-1}   \leq  y_1 \leq 2(X_1A_1)^{-1}  \quad (X_2 A_2)^{-1}  \leq y_2 \leq 2(X_2A_2)^{-1}, \\
  & (2X_1X_2A_2^2)^{-1}  \leq \xi_1^{1/2} \xi_2^{-1} \leq 2(X_1X_2A_2^2)^{-1} , \quad (2X_1X_2A_1^{2})^{-1}   \leq \xi_2^{1/2} \xi_1^{-1} \leq 2(X_1X_2A_1^2)^{-1}.  
   \end{split} 
\end{equation} 
This implies   
\begin{equation}\label{Xi}
  \xi_2 \asymp \Xi_2 := A_1^{4/3}A_2^{8/3} (X_1X_2)^2, \quad  \xi_1 \asymp \Xi_1 := A_2^{4/3}A_1^{8/3}(X_1X_2)^2
\end{equation}  
  which yields in particular the second part of \eqref{J2} as well as 
\begin{equation}\label{xrange}
\begin{split}
 & x_1 \ll A_2^{2/3}A_1^{4/3}X_1X_2, \quad x_2  \ll A_1^{2/3}A_2^{4/3}X_1X_2, \quad x_3\ll A_1^{2/3}A_2^{4/3}X_1X_2, \\
  & x_1x_2-x_3  \ll A_2^{2/3}A_1^{4/3}X_1X_2. 
\end{split}  
\end{equation}
For future purposes we study the volume of the set of $(x_1, x_2, x_3)$ defined by \eqref{Xi} or by
\begin{equation}\label{xi1}
  \xi_2 =  \Xi_2 (1+O(1/R_2)), \quad \xi_1 = \Xi_1(1+O(1/R_1)). 
  \end{equation}  
\begin{lemma}\label{lem4}
 For $\Xi_1, \Xi_2 \geq 1$ and any $\varepsilon > 0$ we have
\begin{displaymath}
  \underset{\substack{x_1, x_2, x_3\\\text{satisfying } \eqref{Xi}}}{\int\int\int} dx_1\, dx_2\, dx_3 \ll (\Xi_1\Xi_2)^{1/2+\varepsilon}.
\end{displaymath}
Moreover,  
\begin{displaymath}  
  \underset{\substack{x_1, x_2, x_3\\\text{satisfying } \eqref{xi1}}}{\int\int\int} dx_1\, dx_2\, dx_3 \ll  (\Xi_1\Xi_2R_1R_2)^{\varepsilon}  \frac{(\Xi_1\Xi_2)^{1/2}}{R_1R_2}. 
\end{displaymath}
\end{lemma}
We postpone the proof to the end of this section. A trivial estimation now implies
\begin{displaymath}
\begin{split}
  \mathcal{J}_{\epsilon_1, \epsilon_2}(A_1, A_2)&  \ll \frac{R_1R_2(R_1+R_2)}{(A_1A_2)^2} \underset{\substack{x_1, x_2, x_3\\\text{satisfying } \eqref{Xi}}}{\int\int\int} dx_1\, dx_2\, dx_3\\
   &\ll R_1R_2(R_1+R_2)   (X_1X_2)^{2}((X_1+A_1)(X_2+A_2))^{\varepsilon}.  
 \end{split} 
\end{displaymath}
Alternatively we write
\begin{equation}\label{exact}
\begin{split}
& \mathcal{J}_{\epsilon_1, \epsilon_2}(A_1, A_2) =  \frac{R_1R_2(R_1+R_2)}{(A_1A_2)^{2} A_2^{  i(\tau_1 - \tau_2)} A_1^{  i(\tau_2-\tau_1)}} \int_0^{\infty}\int_0^{\infty}  (y_1y_2)^{-3i(\tau_1+\tau_2)} \\
  &\times \int_{-\infty}^{\infty}\int_{-\infty}^{\infty}\int_{-\infty}^{\infty} \xi_1^{ -\frac{3}{2}i \tau_1} \xi_2^{- \frac{3}{2} i \tau_2}  e\left(-\epsilon_1 A_1x_1y_1 - \epsilon_2 A_2x_2y_2\right) e\left(-\frac{A_2}{y_2} \cdot \frac{x_1x_3+x_2}{x_3^2 + x_2^2 + 1}\right)\\
 & \times   e\left(-\frac{A_1}{y_1} \cdot \frac{x_2(x_1x_2-x_3) + x_1} {(x_1x_2-x_3)^2 + x_1^2 + 1}\right)\overline{f}(X_1A_1y_1)f\left(\frac{X_1A_2}{y_2} \cdot \frac{\sqrt{(x_1x_2-x_3)^2 + x_1^2 + 1}}{x_3^2 + x_2^2 + 1}\right)\\
& \times \overline{f}(X_2 A_2 y_2)   f\left( \frac{X_2A_1}{y_1} \cdot \frac{ \sqrt{x_3^2 + x_2^2 + 1}} {(x_1x_2-x_3)^2 + x_1^2 + 1}\right)
 dx_1\, dx_2\, dx_3 \frac{dy_1\, dy_2}{y_1y_2}
 \end{split}
\end{equation}
using the notation \eqref{xi}. We can assume that $C_1, C_2$ are integers. Integrating by parts $C_1$ times with respect to $y_1$ and $C_2$ times with respect to $y_2$ introduces an additional factor
\begin{displaymath}
\begin{split}
 \ll_{C_1, C_2}& \left(\left(\frac{y_1}{\tau_1+\tau_2}\right)\left(A_1|x_1| + \frac{A_1|x_2(x_1x_2-x_3)+x_1|}{y_1^2\xi_1} + \frac{1}{y_1}\right)\right)^{C_1} \\
 & \times \left(\left(\frac{y_2}{\tau_1+\tau_2}\right)\left(A_2|x_2| + \frac{A_2|x_1x_3 +x_2|}{y_2^2\xi_2} + \frac{1}{y_2}\right)\right)^{C_2}
 \end{split}
  \end{displaymath}
 and \eqref{J2} follows by \eqref{xrange}  and the same argument that led to \eqref{J1}. \\

The proof of \eqref{improvedJ2} is a small variant of the preceding argument. We need to save an additional power of $R_1+R_2$ which comes from a more careful treatment of the $y_1, y_2$-integral. Let $\eta > 0$ be small. If $1+A_1^{4/3} A_2^{2/3}X \leq (R_1+R_2)^{1-\eta}$, we replace $C_1$ in \eqref{J2} by $C_1 + \eta^{-1}$ saving a factor 
\begin{displaymath}
  \left(\frac{1+A_1^{4/3}A_2^{2/3}X}{R_1+R_2}\right)^{-1/\eta} \geq R_1+R_2.
\end{displaymath}
The same argument works if $1+A_1^{2/3} A_2^{4/3}X \leq (R_1+R_2)^{1-\eta}$. In the remaining case 
\begin{displaymath}
  1+ A_1^{2/3}A_2^{4/3}X \geq (R_1+R_2)^{1-\eta}\quad \text{and} \quad 1+A_2^{2/3}A_1^{4/3}X \geq (R_1+R_2)^{1-\eta}
\end{displaymath}  
   it is enough to show 
   \begin{equation}\label{enough1}
\begin{split}
  \mathcal{J}_{\epsilon_1, \epsilon_2}(A_1,  A_2) d\tau_1 d\tau_2  \ll R_1R_2X^2 ((R_1+R_2)X)^{\varepsilon}; 
  \end{split}
\end{equation}
then the bound \eqref{improvedJ2} follows  with $\varepsilon + \eta(C_1+C_2)$ instead of $\varepsilon$.  To this end, we combine as before \eqref{exact} and the first part of Lemma \ref{lem4}, and need to show that the $y_1$ and $y_2$ integral in \eqref{exact} are both $\ll (\tau_1+\tau_2)^{-1/2} \asymp (R_1+R_2)^{-1/2}$. Our present assumption $X_1 = 1$, $X_2 = X \geq 1$, $A_1, A_2 \leq 1$ together with the size constraints \eqref{sizey} -- \eqref{xrange} imply that the $y_1$ integral is of the form
\begin{equation}\label{statphase}
  \int_0^{\infty} y_1^{-3i(\tau_1+\tau_2)} e(-\epsilon_1A_1x_1y_1) w(A_1y_1) \frac{dy_1}{y_1}
\end{equation}
where $w$ is a smooth function with support in $[1, 2]$ and and $w^{(j)}(y) \ll_j 1$ uniformly in all other variables. We can assume that $\epsilon_1=  \text{sgn}(x_1)$ and $ |x_1| \asymp \tau_1+\tau_2$, otherwise we can save as many powers of $\tau_1+\tau_2$ as we wish by repeated partial integration. In that case we make another change of variables and re-write \eqref{statphase} as
\begin{displaymath}
  \int_0^{\infty} e\left(\frac{3}{2\pi}(\tau_1+\tau_2)(y_1 - \log y_1)\right)    w\left(\frac{3(\tau_1+\tau_2)}{2\pi |x_1|} y_1\right) \frac{dy_1}{y_1}.
\end{displaymath}
A standard stationary phase argument bounds this integral by $(\tau_1+\tau_2)^{-1/2}$:  we cut out smoothly the region $y_1 = 1 + O((\tau_1+\tau_2)^{-1/2})$ which we estimate trivially. For the rest we apply integration by parts. The treatment of the $y_2$ integral is very similar. Here our assumptions imply that the integral is of the form
\begin{displaymath}
\int_0^{\infty}  y_2^{-3i(\tau_1+\tau_2)} 
     e\left(-\frac{A_2}{y_2} \cdot \frac{x_1x_3+x_2}{x_3^2 + x_2^2 + 1}\right)\tilde{w}(XA_2w) \frac{ dy_2}{y_2},
\end{displaymath}
and the same stationary phase-type argument gives a saving of $(\tau_1+\tau_2)^{-1/2}$.\\

Finally we prove \eqref{intJ2}. Let  
\begin{displaymath}
  Z := R_1R_2(X_1+A_1)(X_2+A_2). 
\end{displaymath}
As before we see that we can assume 
 \begin{equation}\label{wlog}
  1+ A_1^{2/3}A_2^{4/3}(X_1+X_2) \geq (R_1+R_2)^{1-\eta}\quad \text{and} \quad 1+A_2^{2/3}A_1^{4/3}(X_1+X_2) \geq (R_1+R_2)^{1-\eta};
\end{equation}  
otherwise we integrate trivially over $\tau_1, \tau_2$. In the situation \eqref{wlog} 
   it is enough to show\footnote{recall that $\mathcal{J}_{\epsilon_1, \epsilon_2}(A_1,  A_2)$ depends on $\tau_1, \tau_2$ although this is not displayed in the notation.}  
\begin{equation}\label{enough}
\begin{split}
\int_{-\infty}^{\infty}\int_{\infty}^{\infty}&g\left(\frac{\tau_1}{R_1}\right)g\left(\frac{\tau_2}{R_2}\right)  \mathcal{J}_{\epsilon_1, \epsilon_2}(A_1,  A_2) d\tau_1 d\tau_2  \ll R_1R_2(R_1+R_2) (X_1X_2)^{2} Z^{\varepsilon}; 
\end{split}
\end{equation}
then the bound \eqref{intJ2} follows  with $\varepsilon + \eta(C_1+C_2)$ instead of $\varepsilon$.   In order to show \eqref{enough},  we integrate \eqref{exact} explicitly over $\tau_1, \tau_2$ and observe that  
\begin{displaymath}
\begin{split}
  &\int_{-\infty}^{\infty}\int_{\infty}^{\infty} g\left(\frac{\tau_1}{R_1}\right)g\left(\frac{\tau_2}{R_2}\right) \frac{ A_1^{  i(\tau_1 - \tau_2)} A_2^{  i(\tau_2-\tau_1)}}{ (y_1y_2)^{3i(\tau_1+\tau_2)}  \xi_1^{ \frac{3}{2}i \tau_1} \xi_2^{\frac{3}{2} i \tau_2}  }d\tau_1\, d\tau_2 \\
  &= R_1R_2 \tilde{g}\left(\frac{R_1}{2\pi} \log\frac{A_1}{A_2(y_1y_2)^3\xi_1^{3/2}}\right)\tilde{g}\left(\frac{R_2}{2\pi} \log\frac{A_2}{A_1(y_1y_2)^3\xi_2^{3/2}}\right). 
 \end{split} 
\end{displaymath}
Since $g$ is smooth, $\tilde{g}$ is rapidly decaying, and up to a negligible error of $Z^{-A}$ we can restrict $\xi_1$, $\xi_2$ to
 \begin{equation}\label{eq}
  \xi_1  = \frac{A_1^{2/3}}{A_2^{2/3} (y_1y_2)^{2}}\left(1 + O\left(\frac{1}{R_1Z^{\varepsilon}}\right)\right), \quad  \xi_2  = \frac{A_2^{2/3}}{A_1^{2/3} (y_1y_2)^{2}}\left(1 + O\left(\frac{1}{R_2Z^{\varepsilon}}\right)\right).
\end{equation}
We note that $ A_1^{2/3}(A_2^{2/3} (y_1y_2)^{2})^{-1} \asymp \Xi_1$, $ A_2^{2/3}(A_1^{2/3} (y_1y_2)^{2})^{-1} \asymp \Xi_2$ in the notation of \eqref{Xi}. Hence a trivial estimate bounds the left hand side of \eqref{enough} by
\begin{displaymath}
  \ll \frac{(R_1R_2)^{2}(R_1+R_2)}{(A_1A_2)^2}   \underset{\substack{x_1, x_2, x_3\\\text{satisfying } \eqref{eq}}}{\int\int\int} dx_1\, dx_2\, dx_3
\end{displaymath}
and the desired bound follows from the second part of Lemma \ref{lem4}. \qed\\

It remains to prove Lemma \ref{lem4}: the conditions $\xi_1 \leq \Xi_1, \, \xi_2 \leq \Xi_2$ are equivalent to
\begin{displaymath}
\begin{split}
 & x_2^2 \leq \Xi_2-1, \quad x_3^2 \leq \Xi_2 - 1 - x_2^2, \quad  \left(x_1 - \frac{x_2x_3}{x_2^2+1}\right)^2 \leq \frac{\Xi_1-1}{x_2^2+1} - \frac{x_3^2}{(x_2^2+1)^2}.    
\end{split}
\end{displaymath}
Hence
\begin{displaymath}
  \underset{\substack{\xi_1 \leq \Xi_1,  \xi_2 \leq \Xi_2}}{\int\int\int} dx_1 \, dx_2 \, dx_3 \leq 8 \int_0^{\Xi_2^{1/2}} \int_0^{\Xi_2^{1/2}} \frac{\Xi_1^{1/2}}{(x_2^2+1)^{1/2}} dx_3\, dx_2 \ll (\Xi_1\Xi_2)^{1/2} \log(1+\Xi_2). 
\end{displaymath}
This proves the first part of the lemma. The second part is more technical. The conditions \eqref{xi1} imply
\begin{equation}\label{imply}
\begin{split}
 & x_2^2 \leq \Xi_2\left(1 +\frac{c}{R_2}\right)-1 , \quad \Xi_2\left(1 - \frac{c}{R_2}\right) -1- x_2^2  \leq x_3^2 \leq \Xi_2\left(1 + \frac{c}{R_2}\right) -1- x_2^2, \\
 & \frac{\Xi_1(1- c/R_1)-1}{x_2^2+1} - \frac{x_3^2}{(x_2^2+1)^2} \leq \left(x_1 - \frac{x_2x_3}{x_2^2+1}\right)^2 \leq \frac{\Xi_1(1+c/R_1)-1}{x_2^2+1} - \frac{x_3^2}{(x_2^2+1)^2}   
\end{split}
\end{equation}
for some constant $c > 0$. We separate four cases for the range of $x_3$. \\

\emph{Case 1.} If $(\Xi_1(1+c/R_1)-1)(x_2^2 + 1)  < x_3^2$, then the condition on $x_1$ is empty. \\

\emph{Case 2.} Let us assume
\begin{equation}\label{x31}
 ( \Xi_1(1-c/R_1)-1)(x_2^2 + 1) \leq x_3^2 \leq (\Xi_1(1+c/R_1)-1)(x_2^2 + 1).
\end{equation}  
Then the volume of the $x_1$-region is 
\begin{displaymath} 
\leq 2 \left(\frac{\Xi_1(1+c/R_1)-1}{x_2^2+1} - \frac{x_3^2}{(x_2^2+1)^2} \right)^{1/2} \leq 2 \left(\frac{2c\, \Xi_1}{R_1(x^2_2+1)}\right)^{1/2}.
\end{displaymath}
The region \eqref{x31} and the second inequality in \eqref{imply} have a non-empty intersection only if
\begin{equation}\label{x21}
  \frac{\Xi_2(1-c/R_2) }{\Xi_1(1+c/R_1)} -1\leq x_2^2 \leq   \frac{\Xi_2(1+c/R_2)  }{\Xi_1(1-c/R_1)}-1. 
\end{equation}
If $R_1, R_2 > 2c$ (which we may assume),  this condition is empty unless $\Xi_2 \geq \frac{1}{3}\Xi_1$. This implies $x_2^2  \leq 3\Xi_2/\Xi_1$ and $x_2^2+1 \asymp \Xi_2/\Xi_1$. In the following we will frequently use the inequality $\sqrt{A}-\sqrt{B} \leq (A-B)A^{-1/2}$ for $A \geq B\geq 0$. Since we are assuming that $\Xi_1, \Xi_2$ are sufficiently large, we can deduce from \eqref{x31} and the second inequality in \eqref{imply}  that the volume of the $x_3$-region is
\begin{equation}\label{before}
  \ll \min\left(\frac{\Xi_2^{1/2}}{R_2}, \frac{\Xi_2^{1/2}}{R_1} \right) 
\end{equation} 
(uniformly in $x_2$) and hence the total contribution under the assumption \eqref{x31} is
\begin{displaymath}
 \int_{x_2 \in \eqref{x21}} \frac{(\Xi_1\Xi_2)^{1/2}}{R_1^{1/2}(R_1+R_2) (x_2^2+1)^{1/2}} dx_2 \asymp  \int_{x_2 \in \eqref{x21}} \frac{\Xi_1}{R_1^{1/2}(R_1+R_2) } dx_2.
\end{displaymath}
The region \eqref{x21} describes an interval of length $O( (R_1^{-1} + R_2^{-1}) \Xi_2/\Xi_1)$ for $x_2^2$, hence the total contribution is
\begin{displaymath}
\ll   \frac{\Xi_1}{R_1^{1/2} (R_1+R_2)}  \frac{\Xi_2^{1/2}}{\Xi_1^{1/2}}\left(\frac{1}{R_1^{1/2}} + \frac{1}{R_2^{1/2}}\right) \ll \frac{(\Xi_1\Xi_2)^{1/2}}{R_1R_2}
\end{displaymath}
as claimed.\\
 
\emph{Case 3.} For a parameter $1/3 \leq \alpha \leq c/R_1$ consider the region
\begin{equation}\label{x32}
  ( \Xi_1(1-2\alpha)-1)(x_2^2 + 1) \leq x_3^2 \leq (\Xi_1(1-\alpha)-1)(x_2^2 + 1).
\end{equation}
The procedure here is very similar to case 2. The $x_1$-volume is at most
\begin{displaymath}
  \leq  \frac{4\Xi_1c/R_1}{(x_2^2+1)^{1/2}} \left(\Xi_1\left(1+\frac{c}{R_1}\right)-1 - \frac{x_3^2}{x^2_2+1}\right)^{-1/2} \ll \frac{\Xi_1^{1/2}}{R_1(x_2^2+1)^{1/2}\alpha^{1/2}}. 
\end{displaymath}
The region \eqref{x32} and the second inequality in \eqref{imply} have a non-empty intersection only if
\begin{equation}\label{x22a}
  \frac{\Xi_2(1-c/R_2)  }{\Xi_1(1-\alpha)} -1 \leq x_2^2 \leq   \frac{\Xi_2(1+c/R_2)  }{\Xi_1(1-2\alpha)} -1. 
\end{equation}
In particular this implies $x_2^2+1 \asymp \Xi_2/\Xi_1$. 
As in \eqref{before} we see that the $x_3$-volume is $\ll \Xi_2^{1/2} \min(R_2^{-1}, \alpha)$, hence the total contribution in the present subcase is
\begin{displaymath}
  \ll \int_{x_2 \in \eqref{x22a}} \frac{\Xi_1 \min(R_2^{-1}, \alpha)}{R_1\alpha^{1/2}} dx_2 \ll  \frac{\Xi_1 \min(R_2^{-1}, \alpha)}{R_1\alpha^{1/2}}    \frac{\Xi_2^{1/2}}{\Xi_1^{1/2}}\left(\frac{1}{R_2^{1/2}} + \alpha^{1/2}\right) \ll \frac{(\Xi_1\Xi_2)^{1/2}}{R_1R_2}.  
  \end{displaymath}\\

\emph{Case 4.} Finally we consider the region $x_3^2 \leq  (\Xi_1/3-1) (x_2^2 + 1)$. In this case  the $x_1$-volume is
\begin{displaymath}
  \leq \frac{4\Xi_1c/R_1}{(x_2^2+1)^{1/2}} \left(\Xi_1\left(1+\frac{c}{R_1}\right)-1 - \frac{x_3^2}{x^2_2+1}\right)^{-1/2} \ll \frac{\Xi_1^{1/2}}{R_1(x_2^2+1)^{1/2}}. 
\end{displaymath}
 The length of the $x_3$ interval  
is at most 
\begin{displaymath}
\leq   \frac{4c\Xi_2}{R_2(\Xi_2(1+c/R_2) -1 - x_2^2)^{1/2}}.
\end{displaymath}  
 Hence the total contribution is at most
\begin{displaymath}
\begin{split}
& \ll  \int_{x_2^2 \leq \Xi_2(1+c/R_2) -1} \frac{\Xi_1^{1/2}\Xi_2\, dx_2}{R_1R_2(x_2+1)^{1/2}  (\Xi_2(1+c/R_j) -1- x_2^2 )^{1/2}}\\
 & \ll \frac{(\Xi_1\Xi_2)^{1/2}}{R_1R_2} +  \int_{1 \leq x_2^2 \leq \Xi_2(1+c/R_2) -1} \frac{\Xi_1^{1/2}\Xi_2\, dx_2}{R_1R_2x (\Xi_2(1+c/R_j) -1- x_2^2 )^{1/2}}.
 \end{split}
\end{displaymath} 
 This last integral can be computed explicitly:
\begin{displaymath}
  \int \frac{dx}{x(Z-x^2)^{1/2}} = \frac{\log(x) - \log(Z + \sqrt{Z(Z-x^2)})}{\sqrt{Z}}, 
\end{displaymath}
and the desired bound follows. This completes the proof of the lemma. \qed\\

 \section{Proofs of the Theorems}
 
For the \textbf{proof of Theorem \ref{thm1}} we choose $n_1 = n_2 = m_1 = m_2=1$ and combine \eqref{spectral}, \eqref{arithmetic}, Lemma \ref{lem1} and Propositions \ref{whitint} and \ref{prop6}.  We choose $\tau_1 = R_1 = T_1$, $\tau_2 = R_2 = T_2$ and  $X_1=X_2 = 1$  in \eqref{defE},  fix a function $f$ and drop all these parameters from the  notation of $F$. By the second part of \eqref{J1} and \eqref{J2}, the Kloosterman terms $\Sigma_{2a}, \Sigma_{2b}, \Sigma_3$ are finite sums over $D_1, D_2$ and hence  are $O((T_1T_2)^{-100})$ by the first part of \eqref{J1} and \eqref{J2}. 
 The diagonal term \eqref{diag} is $\asymp T_1T_2(T_1+T_2)$. On the spectral side, we drop the Eisenstein spectrum and large parts of the cuspidal spectrum to conclude by \eqref{asymp} and Lemma \ref{lem1} the upper bound
\begin{displaymath}
   \sum_{\substack{| \nu^{(j)}_1 -iT_1| \leq c\\  | \nu^{(j)}_2 - iT_2| \leq c}} \left(\underset{s=1}{\rm res} L(s, \phi_j \times \tilde{\phi}_j)\right)^{-1}  \ll T_1T_2(T_1+T_2)
\end{displaymath}
for some sufficiently small $c$ and $T_1, T_2 \geq T_0$, and hence
\begin{equation}\label{prelimupper}
   \sum_{\substack{| \nu^{(j)}_1 -iT_1| \leq K\\  | \nu^{(j)}_2 - iT_2| \leq K}} \left(\underset{s=1}{\rm res} L(s, \phi_j \times \tilde{\phi}_j)\right)^{-1}  \ll_K T_1T_2(T_1+T_2)
\end{equation}
for any $K \geq 1$ by adding the contribution of $O_K(1)$ balls.  To prove the lower bound, we   choose (once and for all) $K$  so large that
\begin{displaymath}
    \sum_{\substack{\max(| \nu^{(j)}_1 -iT_1|,  | \nu^{(j)}_2 - iT_2|)  \geq K}}  \|\phi_j\|^{-2} |\langle \tilde{W}_{\nu_1, \nu_2}, F\rangle |^2  \leq \frac{1}{2} \| F \|^2
  \end{displaymath}
  which is possible by \eqref{prelimupper} and \eqref{bound}. We bound the Eisenstein spectrum trivially:  the second line of \eqref{spectral} contributes $O((T_1+T_2)^{\varepsilon})$ by known bounds for the zeta function on the line $\Re s = 1$, the third line contributes similarly  $O((T_1+T_2)^{1+\varepsilon})$ by Weyl's law for $SL_2(\Bbb{Z})$ and lower bounds for the $L$-functions in  the denomiator \cite{HL, HR}. Hence we obtain
\begin{displaymath}
    \sum_{\substack{ | \nu^{(j)}_1 -iT_1| \leq K\\   | \nu^{(j)}_2 - iT_2|  \leq K}}  \|\phi_j\|^{-2} |\langle \tilde{W}_{\nu_1, \nu_2}, F\rangle |^2  \geq \frac{1}{2} \| F \|^2 + O((T_1+T_2)^{1+\varepsilon}),
  \end{displaymath}  
 and the lower bound in Theorem \ref{thm1} follows from \eqref{bound} and  \eqref{diag} for $T_1, T_2$ sufficiently large. \qed\\
  
The \textbf{proof of Theorem \ref{thm2}} proceeds similarly. As mentioned in the introduction, as a direct corollary of Theorem \ref{thm1} we find that the number of exceptional Maa{\ss} forms $\phi_j$ with $\gamma_j = T+O(1)$ is $O(T^2)$. In order to prove Theorem \ref{thm2}, it is therefore enough to consider those Maa{\ss} forms with $|\rho_j| \geq \varepsilon$. Moreover, by symmetry it is enough to bound only Maa{\ss} forms satisfying \eqref{nontemp}. In \eqref{defE} we take $\tau_2 = R_2 = T$, $R_1 = 1$, $\tau_1 = 0$, $X_1 = 1$, $X_2 = X = T^{\delta}$ for some $\delta > 0$ to be chosen later. With this data, the spectral side, after dropping 
\begin{itemize}
\item   the tempered spectrum, 
\item  the Eisenstein spectrum and
\item  those parts of the non-tempered spectrum not of the form \eqref{nontemp} with $|\rho_j| \geq \varepsilon$, 
\end{itemize}
is  by \eqref{asymp1} (note that \eqref{nontempasump} is satisfied) and the upper bound of \eqref{resiweak}  at least
\begin{displaymath}
\gg  T^{-\varepsilon} X^2    \sum_{\substack{ \phi_j \text{ as in } \eqref{nontemp}\\ \gamma_j = T+O(1)\\ |\rho_j| \geq \varepsilon}}    X^{2|\rho_j|}. 
\end{displaymath} 
On the arithmetic side, the diagonal term is $\asymp T^2X^2$ by \eqref{diag}. Next by \eqref{J1} we have
\begin{displaymath}
  \Sigma_{2a} \ll T^2X  \sum_{D_1 \ll 1} \frac{|\tilde{S}(\pm 1, 1, 1, D_1, D_1^2)|}{D_1^3} T^{-102} \ll X T^{-100}
\end{displaymath}
and
\begin{displaymath}
  \Sigma_{2b} \ll T^2X^2 \sum_{D_2 \ll X^{1/2}}  \frac{|\tilde{S}(\pm 1, 1, 1, D_2, D_2^2)|}{D_2^3} T^{-102} \ll X^{2+\varepsilon} T^{-100}
\end{displaymath}
by \eqref{sharp}. (Note that we are exchanging $X_1$ and $X_2$ for $\Sigma_{2b}$.)  
The long Weyl element contributes at most
\begin{equation}\label{longWeyl}
 (TX^2)^{1+\varepsilon}  \sum_{D_1, D_2 \ll X}  \frac{|S(\pm 1, 1, 1, 1, D_1, D_2)|}{D_1 D_2}  \left(\frac{1+X/D_2}{T}\right)^{-C_1}\left(\frac{1+X/D_1}{T}\right)^{-C_2} \ll (TX^2)^{1+2\varepsilon}  \frac{X}{T}
\end{equation}
which follows by combining \eqref{improvedJ2} and \eqref{kloosum}.  Choosing $X=T^2$ completes the proof of Theorem \ref{thm2}. \qed\\ 

We proceed to prove \textbf{Theorem \ref{thm3}}. Again we choose $X_1 = X_2 = 1$, $R_1 = T_1$, $R_2 = T_2$ in \eqref{defE}, fix a function $f$ and then drop $R_1, R_2, X_1, X_2, f$ from the notation of $F$ and keep only $\tau_1, \tau_2$.  We also fix a suitable non-negative smooth function $g$ with support in $[1/2, 3]$ as in Proposition \ref{prop6}.  Let $T := \max(T_1, T_2)$. The left hand side of \eqref{large} is, by \eqref{asymp} and the upper bound of \eqref{resiweak}, 
\begin{displaymath}
 \ll T^{\varepsilon}  \sum_{j}  \frac{1}{\|\phi_j\|^2} \int_{T_1}^{2T_1} \int_{T_2}^{2T_2}  |\langle \tilde{W}_{\nu_1, \nu_2}, F_{\tau_1, \tau_2}\rangle|^2 d\tau_1 d\tau_2  \Bigl| \sum_{n \leq N} \alpha(n) A_j(n, 1)\Bigr|^2.
\end{displaymath}
We cut the $n$-sum   into   dyadic intervals, insert artificially the function $g$  and bound the preceding display by
\begin{displaymath}
 \ll (NT)^{\varepsilon} \max_{\substack{M \leq N }} \sum_{j}  \frac{1}{\|\phi_j\|^2} \int_{-\infty}^{-\infty} \int_{-\infty}^{\infty}  g\left(\frac{\tau_1}{T_1}\right)g\left(\frac{\tau_2}{T_2}\right)    |\langle \tilde{W}_{\nu_1, \nu_2}, F_{\tau_1, \tau_2}\rangle|^2 d\tau_1 d\tau_2   \Bigl| \sum_{M \leq n \leq 2M} \alpha(n) A_j(n, 1)\Bigr|^2.
\end{displaymath}
Next we add artificially the continuous spectrum getting the upper bound
\begin{displaymath}
\begin{split}
&(NT)^{\varepsilon} \max_{ \substack{M \leq N}} \left(\sum_{j}  \int_{-\infty}^{\infty} \int_{-\infty}^{\infty} g\left(\frac{\tau_1}{T_1}\right)g\left(\frac{\tau_2}{T_2}\right)  \frac{ |\langle \tilde{W}_{\nu_1, \nu_2}, F_{ \tau_1, \tau_2}\rangle|^2}{\|\phi_j\|^2} d\tau_1 d\tau_2    \Bigl| \sum_{M \leq n \leq 2M} \alpha(n) A_j(n, 1)\Bigr|^2\right.\\
&+  \frac{1}{(4\pi i)^2} \int_{(0)} \int_{(0)}  \int_{-\infty}^{\infty} \int_{-\infty}^{\infty}   \frac{g\left(\frac{\tau_1}{T_1}\right)g\left(\frac{\tau_2}{T_2}\right)|\langle \tilde{W}_{\nu_1, \nu_2}, F_{ \tau_1, \tau_2}\rangle|^2}{|\zeta(1+3\nu_0)\zeta(1+3\nu_1)\zeta(1+3\nu_2)|^2} d\tau_1 d\tau_2  
 \Bigl| \sum_{M \leq n \leq 2M} \alpha(n) A_{\nu_1, \nu_2}(n, 1)\Bigr|^2 d\nu_1 d\nu_2\\
 & \left.+    \frac{c}{2\pi i} \sum_{j} \int_{(0)}  \int_{-\infty}^{\infty} \int_{-\infty}^{\infty}   \frac{ g\left(\frac{\tau_1}{T_1}\right)g\left(\frac{\tau_2}{T_2}\right)|\langle \tilde{W}_{\frac{2}{3}\nu_j, \mu - \frac{1}{3} \nu_j}, F_{ \tau_1, \tau_2}\rangle|^2 }{|L(1+3\mu, u_j)|^2 L(1, \text{Ad}^2 u_j)\ }   d\tau_1 d\tau_2  \Bigl| \sum_{M \leq n \leq 2M} \alpha(n) B_{\mu, u_j}(n, 1)\Bigr|^2
  d\mu\right). 
    \end{split}
\end{displaymath}
We open the squares and apply the Kuznetsov formula, that is, we replace the three terms of the shape \eqref{spectral} with the four terms \eqref{arithmetic}. We estimate each of them individually. The diagonal term contributes by \eqref{diag}
\begin{displaymath}
 \ll  (NT)^{\varepsilon} \max_{ \substack{M \leq N }} \sum_{M \leq n \leq 2M} |\alpha(m)|^2 \int_{-\infty}^{\infty} \int_{-\infty}^{\infty}  g\left(\frac{\tau_1}{T_1}\right)g\left(\frac{\tau_2}{T_2}\right) T_1T_2(T_1+T_2) d\tau_1 d\tau_2 \ll (NT)^{\varepsilon} T_1^2T_2^2 (T_1+T_2)  \| \alpha \|_2^2.  
\end{displaymath}
This is the first term on the right hand side of \eqref{large}. In the term $\Sigma_{2a}$ in \eqref{arithmetic} the condition $D_1 \mid D_2$, $D_1^2 = nD_2$ is equivalent to $D_1 = nd$, $D_2 = nd^2$ for some $d \in \Bbb{N}$; hence its contribution is at most  
\begin{displaymath}
\begin{split}
& \ll  (NT)^{\varepsilon}  \max_{\substack{ M \leq N}} \sum_{M \leq n, m \leq 2M}| \alpha(n)  \alpha(m)|\sum_{\epsilon = \pm 1} \sum_{d} \frac{|S(\epsilon m,n, 1, nd, nd^2)|}{n^2d^3} \\
 &\quad\quad \times  \left| \int_{-\infty}^{\infty}\int_{-\infty}^{\infty}  g\left(\frac{\tau_1}{T_1}\right)g\left(\frac{\tau_2}{T_2}\right)\tilde{\mathcal{J}}_{\epsilon, F}\left( \frac{m^{1/2}}{nd^{3/2}} \right)d\tau_1 d\tau_2 \right|.\\
    \end{split}
\end{displaymath}
By \eqref{J1}, the $d$-sum is finite, hence in combination with \eqref{sharp} this is bounded by
\begin{displaymath}
 \ll (NT)^{\varepsilon}  \max_{\substack{ M \leq N}} \sum_{M \leq n, m \leq 2M}| \alpha(n)  \alpha(m)| T^{-101} \ll N^{\varepsilon} T^{-100} \|\alpha\|^2.
\end{displaymath}
 
 In the term $\Sigma_{2b}$ in \eqref{arithmetic} the condition $D_2 \mid D_1$ is redundant, and the argument of $\tilde{\mathcal{J}}_{\epsilon, F^{\ast}}$ equals $(n/(mD_2^3))^{1/2}$. As before we see that this contributes at most $N^{\varepsilon} T^{-100} \|\alpha\|^2$.

Finally the long Weyl element finally contributes by \eqref{intJ2}
\begin{equation}\label{long}
\begin{split}
 &\ll  (NT)^{\varepsilon}  \max_{\substack{ M \leq N }} \sum_{M \leq n, m \leq 2M} |\alpha(n)  \alpha(m)|\sum_{\epsilon_1, \epsilon_2 = \pm 1} \sum_{D_1, D_2} \frac{|S(\epsilon_1 m, \epsilon_2 ,n, 1, D_1, D_2)|}{D_1D_2} \\
 & \quad\quad\quad\quad \times  \left| \int_{-\infty}^{\infty}\int_{-\infty}^{\infty}  g\left(\frac{\tau_1}{T_1}\right)g\left(\frac{\tau_2}{T_2}\right)\mathcal{J}_{\epsilon_1, \epsilon_2}\left(\frac{\sqrt{mD_1}}{D_2}, \frac{\sqrt{nD_2}}{D_1}  \right)d\tau_1 d\tau_2 \right|\\
 &\ll_{C_1, C_2} (NT)^{\varepsilon} \max_{\substack{ M \leq N}}  T_1T_2(T_1+T_2) \sum_{M \leq n, m\leq 2M} |\alpha(n)  \alpha(m)|\\
 & \quad\quad\quad\quad \times \sum_{\epsilon_1, \epsilon_2 = \pm 1} \sum_{D_1, D_2} \frac{|S(\epsilon_1 m, \epsilon_2 ,n, 1, D_1, D_2)|}{D_1D_2}     \left(\frac{1 + M/D_2}{T_1+T_2}\right)^{C_1} \left(\frac{1 + M/D_1}{T_1+T_2 }\right)^{C_2}  
    \end{split}
\end{equation}
for any $C_1, C_2 \geq 0$. Recalling the notation $T= \max(T_1, T_2)$  and using \eqref{kloosum}, it is straightforward to see that the pervious display is 
\begin{displaymath}
\begin{split}
 & \ll (NT)^{\varepsilon} \left(T_1T_2N^2 \right) \| \alpha \|_2^2.
\end{split}
\end{displaymath}
 This is the second term on the right hand side of \eqref{large}.   \qed \\

Finally we prove \textbf{Theorem \ref{thm4}}. To this end, we express $L(\phi_j, 1/2)$ by an  approximate functional equation. As we are summing over the archimedean parameters of the $L$-functions, we need an approximate functional equation whose weight function is essentially independent of the underlying family. This has been obtained in \cite[Proposition 1]{BH}, and we quote the following special case. For a Maa{\ss} form $\phi_j$ as in Theorem \ref{thm4} put 
\begin{displaymath}
\begin{split}
 & (\eta_j)_1 = \frac{1}{4} + \frac{2\nu^{(j)}_1 +\nu^{(j)}_2}{2}, \quad (\eta_j)_2 = \frac{1}{4} + \frac{-\nu^{(j)}_1 + \nu^{(j)}_2}{2}, \quad (\eta_j)_3 = \frac{1}{4} + \frac{-\nu^{(j)}_1 -2 \nu^{(j)}_2}{2}\\
 & \eta_j= \min_{1 \leq l \leq 3} |(\eta_j)_l| \asymp |(\eta_j)_2|, \quad C_j = \prod_{l = 1}^3 |(\eta_j)_l|.
 \end{split}
\end{displaymath}
Moreover, for a multi-index $\textbf{n} \in \Bbb{N}_0^6$ we write $|\textbf{n}| = \textbf{n}(1) + \ldots + \textbf{n}(6)$ and
\begin{displaymath}
  {\bm \eta}_j^{-\textbf{n}} := \prod_{l=1}^3 (\eta_j)_l^{-\textbf{n}(2l-1)} (\overline{\eta}_j)_l^{-\textbf{n}(2l)}. 
\end{displaymath}

\begin{prop} Let $G_0 : (0, \infty) \rightarrow \RR$ be a smooth function with functional equation $G_0(x) + G_0(1/x) = 1$ and derivatives decaying faster than any negative power of $x$ as $x \rightarrow \infty$. Let $M \in \NN$ and fix a Maa{\ss} form $\phi$ as above. There are explicitly computable rational constants $c_{{\bm n},\ell} \in \QQ$ depending only on ${\bm n}$, $\ell$, $M$ such that the following holds for
\begin{displaymath}
  G(x)  :=   G_0(x)+\sum_{\substack{0< |{\bm n}| < M \\ 0< \ell < |{\bm n}|+M}} c_{{\bm n},\ell} {\bm\eta}_j^{-\bm n}   \left(x\frac{\partial}{\partial x}\right)^{\ell} G_0(x).
\end{displaymath}
For any $\eps > 0$ one has
\begin{equation}\label{approxfuncteq}
  L(\phi_j, 1/2) = \sum_{n=1}^{\infty} \frac{A_j(1,  n)}{\sqrt{n}} G\left(\frac{n}{\sqrt{C_j}}\right) + \kappa_j  \ov{\sum_{n=1}^{\infty} \frac{A_j(1, n)}{\sqrt{n}} G\left(\frac{n}{\sqrt{C_j}}\right)} + O\bigl(\eta_j^{-M}C_j^{1/4+\eps}\bigr),
\end{equation}
where $|\kappa_j| = 1$ and the implied constant depends at most on $\eps$, $M$,  and the function $G_0$.  
\end{prop}
 
It is now a simple matter to prove Theorem \ref{thm4}.  We can assume that $T$ is sufficiently large.  Let 
\begin{displaymath}
  G_{\ell}(x) :=  \left(x\frac{\partial}{\partial x}\right)^{\ell} G_0(x).
\end{displaymath}
Then the Mellin transform $\widehat{G}_{\ell}(s)$ is rapidly decaying on vertical lines $\Re s = \sigma > 0$. By \eqref{approxfuncteq} and Mellin inversion we have
\begin{displaymath}
\begin{split}
  & \sum_{\substack{T \leq |\nu^{(j)}_1|, |\nu^{(j)}_2| \leq  2T}} |L(\phi_j, 1/2)|^2 \\
   &\ll_{M, \varepsilon} \sum_{\ell \leq 2M}  \sum_{\substack{T  \leq |\nu^{(j)}_1|,   |\nu^{(j)}_2| \leq  2T}} \Bigl( \Bigl| \sum_{n \leq T^{3/2+\varepsilon}} \frac{A_j(n, 1)}{\sqrt{n}} \overline{G}_{\ell}\Bigl(\frac{n}{\sqrt{C_j}}\Bigr)\Bigr| + O(\eta_j^{-M} C_j^{1/4+\varepsilon} + T^{-100})\Bigr)^2\\
 & \ll_{M, \varepsilon} T^{\varepsilon}      \sum_{\substack{T \leq |\nu^{(j)}_1|,  |\nu^{(j)}_2| \leq  2T} } \Bigl( \int_{-T^{\varepsilon}}^{T^{\varepsilon}} \Bigl| \sum_{n \leq T^{3/2+\varepsilon}} \frac{A_j(n, 1)}{n^{1/2 + \varepsilon + it}}  \Bigr|dt + O(\eta_j^{-M+1/4} T^{1/2+\varepsilon} + T^{-100})\Bigr)^2,
\end{split}   
  \end{displaymath}
noting that $C_j \ll \eta_jT^2 \ll T^3$.  This is at most
  \begin{displaymath}
    \ll T^{\varepsilon} \Bigl(\max_{|t| \leq T^{\varepsilon}}  \sum_{\substack{T  \leq |\nu^{(j)}_1|,  |\nu^{(j)}_2| \leq  2T} }\Bigl| \sum_{n \leq T^{3/2+\varepsilon}} \frac{A_j(n, 1)}{n^{1/2 + \varepsilon +  it}}  \Bigr|^2 +  T \sum_{T \leq |\nu^{(j)}_1|, |\nu^{(j)}_2| \leq  2T} (1 + |\nu^{(j)}_1 - \nu^{(j)}_2|)^{-2M + 1/2} \Bigr).
  \end{displaymath}
By Theorem \ref{thm3} and \eqref{ranksel} the first term is $O(T^{5+\varepsilon})$. By Theorem \ref{thm1} or \eqref{prelimupper}  it is easy to see that the second term is also $O(T^{5  + \varepsilon})$. This completes the proof of Theorem \ref{thm4}. \qed\\

\section{Appendix: A theorem of Goldfeld-Kontorovich}

A very nice  application of the $GL(3)$ Kuznetsov formula has been given recently in \cite{GK2}. The purpose of this appendix is to  illustrate how the methods of this paper directly yield a version of \cite[Theorem 1.3]{GK2} with considerably better error terms and without assuming the Ramanujan conjectures. 
 We keep the notation developed so far. 

\begin{theorem}\label{appendix} For $n_1, n_2, m_1, m_2 \in \Bbb{N}$, $P = n_1n_2m_1m_2 $, $T$ sufficiently large,    one has
\begin{displaymath}
   \begin{split} 
 & \sum_{j}  \overline{A_j(n_1, n_2)} A_j(m_1, m_2) \frac{h_{T}(\nu_1, \nu_2)}{\| \phi_j \|^2} = \delta_{\substack{m_1=n_1\\ m_2=n_2}} \sum_{j}   \frac{h_{T }(\nu_1, \nu_2) }{\|\phi_j \|^2}  + O((T^2P^{1/2} + T^{3}P^{\theta} +P^{5/3})(TP)^{\varepsilon})
  \end{split}
\end{displaymath}
where $\theta \leq 7/64$ is a bound towards the Ramanujan conjecture on $GL(2)$. Here $h_T$ is non-negative, uniformly bounded on $\{|\Re \nu_1| \leq 1/2\} \times \{|\Re \nu_2| \leq 1/2\}$, $h_T \asymp 1$ on $\{(\nu_1, \nu_2) \mid c \leq \Im \nu_1, \Im \nu_2 \leq T, \, |\Re \nu_1|, |\Re \nu_2| \leq 1/2\}$ for some absolute constant $c > 0$, and $h_T(\nu_1, \nu_2) \ll_A ((1+|\nu_1|/T)(1+|\nu_2|/T))^{-A}$. 
 \end{theorem}
 
For comparison, the error term in \cite[Theorem 1.3]{GK2} (scaled down by $T^{-3R}$) is   $O(T^{3+\varepsilon} P^2)$, but see also \cite[Remarks 1.8, 1.19]{GK2} where possible improvements are mentioned.  
A more precise discussion on the asymptotic behaviour of the test function $h_T$ can be found in Remark 6 below. \\ 

Injecting Theorem \ref{appendix} into the estimates of \cite[Section 9]{GK2} and using only $\theta   < 1/3$, we obtain the following corollary. For a Hecke-Maa{\ss} form $\phi$ for $SL_3(\Bbb{Z})$ let $\rho(\phi)$ be one of $\phi$, $\text{sym}^2 \phi$ or $\text{Ad}\phi$. Let $\psi$ be a smooth test function whose Fourier transform has support in $(-\delta, \delta)$ for some $\delta > 0$. Define $D(\rho(\phi), \psi)$ as in \cite[Section 1.4]{GK2}.  
\begin{cor} [Goldfeld-Kontorovich] Assume the generalized Riemann hypothesis and the  Ramanujan conjectures. Suppose 
\begin{displaymath}
\begin{split}
&  \delta < 5/23,\,\,\, \quad \rho(\phi) = {\rm sym}^2\phi \text{ or } {\rm Ad}\phi,\\
 & \delta < 10/13,   \quad \rho(\phi) = \phi,
\end{split}  
\end{displaymath}
Then one has
\begin{displaymath}
 \Bigl( \sum_{j}   \frac{h_{T }(\nu_1, \nu_2) }{\|\phi_j \|^2} \Bigr)^{-1} \sum_j D(\rho(\phi_j), \psi)   \frac{h_{T }(\nu_1, \nu_2) }{\|\phi_j \|^2} = \int_{\Bbb{R}} \psi(x) W_{\rho}(x) dx + O\left(\frac{\log\log T}{\log T}\right),
\end{displaymath}
where 
\begin{displaymath}
  W_{\rho}(x) = 1, \quad \rho(\phi) = \phi \text{ or } \text{sym}^2\phi, \quad W_{\rho}(x) = 1-\frac{\sin(2\pi x)}{2\pi x}, \quad \rho(\phi) = {\rm Ad} \phi.
\end{displaymath}
In particular, the symmetry types are unitary or symplectic, respectively.
\end{cor}

This improves the range of the support of $\hat{\psi}$ by about a factor 3 compared  to \cite[Theorem 1.13]{GK2} (see also \cite[Remarks 1.18, 1.19]{GK2}).\\
 
\textbf{Proof of Theorem \ref{appendix}.}   Let $g$ be a fixed, smooth, non-negative, compactly supported   test function.  Let $R_1,  R_2$ be sufficiently large, and write $R = R_1+R_2$.  We choose $F$ as in \eqref{defE} 
with $X_1 = X_2 = 1$ and integrate the equality in Proposition \ref{kuznetsov} against 
\begin{displaymath}
  \int_0^{\infty} \int_0^{\infty} g\left(\frac{\tau_1}{R_1}\right) g\left(\frac{\tau_2}{R_2}\right) d\tau_1 d\tau_2
\end{displaymath}
as in \eqref{intJ2}. From Proposition \ref{whitint}, the above mentioned  lower bounds for $L$-functions \cite{HL, HR} on the line $\Re s = 1$ and  Weyl's law for $GL(2)$ we conclude that the Eisenstein contribution in \eqref{spectral} is $O(R^{3+\varepsilon} P^{\theta + \varepsilon})$. From Proposition \ref{prop4} and Proposition \ref{prop6} we conclude by the same calculation as in \eqref{long} that the long Kloosterman sum $\Sigma_3$ in \eqref{arithmetic} contributes $O(R^{2+\varepsilon} P^{1/2+\varepsilon}).$ Similarly, if $P < R^{3-\varepsilon}$, the other two Kloosterman contributions $\Sigma_{2a} + \Sigma_{2b}$ are  $O( R^{-100})$, and are otherwise $O(R^{5+\varepsilon})$ which follows after a straightforward estimate using Proposition \ref{prop6} and  \eqref{sharp}.  Hence in either case their contribution is $O(P^{5/3+\varepsilon}).$
 We conclude
\begin{equation}\label{notperfect}
\begin{split}
 & \sum_{j}  \overline{A_j(n_1, n_2)} A_j(m_1, m_2) \frac{h_{R_1, R_2}(\nu_1, \nu_2)}{\| \phi_j \|^2} = \delta_{\substack{m_1=n_1\\ m_2=n_2}} H_{R_1, R_2} + O((R^2P^{1/2} + R^3P^{\theta} + P^{5/3} )(RP)^{\varepsilon})
  \end{split}
\end{equation}
where 
\begin{displaymath}
  h_{R_1, R_2}(\nu_1, \nu_2) =    \int_0^{\infty} \int_0^{\infty} g\left(\frac{\tau_1}{R_1}\right) g\left(\frac{\tau_2}{R_2}\right) |\langle \tilde{W}_{\nu_1, \nu_2}, F\rangle|^2 d\tau_1 d\tau_2
\end{displaymath}
and $H_{R_1, R_2} = \int_0^{\infty} \int_0^{\infty} g (\frac{\tau_1}{R_1} ) g (\frac{\tau_2}{R_2} ) \| F \|^2 d\tau_1 d\tau_2$, but we only need to know that this quantity is independent of $n_1, n_2, m_1, m_2$. 

The weight function $h_{R_1, R_2}$ is uniformly bounded and non-negative. It follows directly from  Proposition \ref{whitint} that 
\begin{equation}\label{weightfunc}
  h_{R_1, R_2}(\nu_1, \nu_2)  \asymp   1 \quad \text{for} \quad |\Re \nu_1|, |\Re \nu_2| \leq 1/2, \quad \frac{\Im \nu_1}{R_1}, \frac{\Im \nu_2}{R_2} \in \text{supp}(g),
\end{equation}
and rapidly decaying outside the region $\frac{|\nu_1|}{R_1}, \frac{|\nu_2|}{R_2} \in \text{supp}(g)$. In other words, $h_{R_1, R_2}$ is a good approximation of the characteristic function on the square $\Im \nu_1 \asymp R_1$, $\Im \nu_2 \asymp R_2$. 

Applying \eqref{notperfect} with $n_1=n_2= m_1= m_2=1$, we see that
\begin{displaymath}
  H_{R_1, R_2} =  \sum_{j}   \frac{h_{R_1, R_2}(\nu_1, \nu_2)}{\| \phi_j \|^2}  + O(R^{3+\varepsilon}). 
\end{displaymath}
Hence we obtain 
\begin{displaymath}
  \begin{split} 
 & \sum_{j}  \overline{A_j(n_1, n_2)} A_j(m_1, m_2) \frac{h_{R_1, R_2}(\nu_1, \nu_2)}{\| \phi_j \|^2} = \delta_{\substack{m_1=n_1\\ m_2=n_2}} \sum_{j}   \frac{h_{R_1, R_2}(\nu_1, \nu_2)}{\| \phi_j \|^2}   + O((R^2P^{1/2} + R^{3}P^{\theta} + P^{5/3}) (RP)^{\varepsilon})
  \end{split}
\end{displaymath}
whenever $R_1, R_2$ are sufficiently large. Piecing together dyadic squares, we obtain Theorem \ref{appendix}. \hfill $\square$\\

\textbf{Remark 6:} The proof of Proposition \ref{whitint} gives much more precise information on the weight function $h_T$ in Theorem \ref{appendix}. By \eqref{asympwhit},  we see that $h_{R_1, R_2}$ described in  \eqref{weightfunc} satisfies the more precise asymptotic relation
\begin{displaymath}
  h_{R_1, R_2}(\nu_1, \nu_2) \sim c \frac{R_1R_2(R_1+R_2)}{|\nu_1\nu_2(\nu_1+\nu_2)|} g\left(\frac{|\nu_1|}{R_1}\right)g\left(\frac{|\nu_2|}{R_2}\right), \quad \nu_1, \nu_2 \in i \Bbb{R},
\end{displaymath}
for $R_1, R_2 \rightarrow \infty$, where the constant $c>0$ is given by
\begin{displaymath}
  c = \frac{(2\pi)^3}{3^3} \int_{\Bbb{R}} \int_{\Bbb{R}} |\widehat{f}(-1-ix-2iy) \widehat{f}(-1-2ix-iy)|^2 dx\, dy = \frac{(2\pi)^3}{3^4} \left( \int_{\Bbb{R}} |\widehat{f}(-1-ix)|^2 dx\right)^2  
  \end{displaymath}
for the weight function $f$ in the Poincar\'e series \eqref{defE}. In particular, by varying $g$ one has the flexibility to prescribe asymptotically any reasonable bump function on the tempered spectrum.

\end{document}